\documentclass[11pt,leqno]{article}
\usepackage{a4}
\usepackage[T1]{fontenc}
\usepackage[english]{babel}
\usepackage{latexsym}
\usepackage{amsmath}
\usepackage{amssymb}
\usepackage{mathrsfs}
\usepackage{color}
\usepackage{cite}
\usepackage{titling}
\makeatletter
\@addtoreset{equation}{section} 
\makeatother
\definecolor{move}{rgb}{.3,.1,.8}
\newcommand{\ud}{\mathrm{d}}
\newcommand{\de}{\mathrm{D}}
\newcommand{\supp}{\mathrm{supp}}
\newcommand{\re}{\mathrm{Re}}

\newcommand{\e}{\mathrm{e}}

\newcommand{\R}{\mathbb{R}}
\newcommand{\N}{\mathbb{N}}

\newcommand{\Ci}{\mathscr{C}^{\infty}}

\newcommand{\fo}{\phi_{1,n}}
\newcommand{\ft}{\phi_{2,n}}
\newcommand{\fot}{\phi_{1/2}}
\newcommand{\fto}{\phi_{2/1}}

\newenvironment{pr}{\vspace{5pt}\textbf{{\small Proof :}}\\}{\hspace{\stretch{1}}\rule{1ex}{1ex}\vspace{5pt}}
\newtheorem{thm}{Theorem}[section]
\newtheorem{pro}{Proposition}[section]

\newtheorem{lem}{Lemma}[section]

\usepackage{fancyhdr}
\fancypagestyle{plain}{\fancyhead{}}
\title{Stability of elastic transmission systems with a local Kelvin-Voigt damping}
\author{{FATHI HASSINE}\\ \textit{UR Analysis and Control of PDE UR13ES64}\\ \textit{Department of Mathematics, Faculty of Sciences of Monastir}\\ \textit{University of Monastir, 5019 Monastir, Tunisia}\\ \textit{email:} \texttt{fathi.hassine@fsm.rnu.tn}}
\date{}
\begin{document}
\maketitle
\begin{center}
\abstract{In this paper, we consider the longitudinal and transversal vibrations of the transmission Euler-Bernoulli beam with Kelvin-Voigt damping distributed locally on any subinterval of the region occupied by the beam and only in one side of the transmission point. We prove that the semigroup associated with the equation for the transversal motion of the beam is exponentially stable, although the semigroup associated with the equation for the longitudinal motion of the beam is polynomially stable. Due to the locally distributed and unbounded nature of the damping, we use a frequency domain method and combine a contradiction argument with the multiplier technique to carry out a special analysis for the resolvent.}
\end{center}


\vspace{10pt}
\section{Introduction}
The theories of viscoelasticity, which include the Maxwell model, the Kelvin-Voigt model, and the standard linear solid model, are used to predict a material's response under different loading conditions. One of the simplest mathematical models constructed to describe the viscoelastic effects is the classical Kelvin-Voigt model. The basic idea concerning this model is that the stress is dependent on the deformation tensor and deformation-rate tensor. This model consists of a Newtonian damper and Hooke's elastic spring connected in parallel.

Recent advances in material science have provided new means for the suppression of vibrations of elastic structures. One approach is to bond or embed patches made of "smart material" to the underlying structure as passive or active damper. Due to the presence of the patches, the material properties of the structure, such as the density, Young's moduli, and damping coefficients, are changed. In particular, jump discontinuities at the location of the edges of the patch are usually introduced into these properties.

Consider a clamped elastic beam of length $L$. One segment of the beam is made of a viscoelastic material with Kelvin-Voigt constitutive relation in which a transmission effect has been established such a way that the damping is locally effective in only one side the transmission boundary. By the Kirchhoff hypothesis and neglecting the rotatory inertia, the longitudinal and transversal vibration of the beam can be described by the following transmission equations and boundary-initial conditions:
\begin{equation}\label{kv1}
\left\{\begin{array}{ll}
\ddot{u}_{1}(x,t)-(p_{1}.u_{1}'+D_{a}\dot{u}_{1}')'(x,t)=0&\text{in }(0,l)\times(0,+\infty)
\\
\ddot{u}_{2}(x,t)-(p_{2}.u_{2}')'(x,t)=0&\text{in }(l,L)\times(0,+\infty)
\\
u_{1}(l,t)=u_{2}(l,t)&\text{for }t\in(0,+\infty)
\\
p_{1}.u_{1}'(l,t)=p_{2}.u_{2}'(l,t)&\text{for }t\in(0,+\infty)
\\
u_{1}(0,t)=u_{2}(L,t)=0&\text{for }t\in(0,+\infty)
\\
u_{1}(x,0)=u_{1}^{0}(x),\;u_{1}'(x,0)=u_{1}^{1}(x)&\text{on }(0,l)
\\
u_{2}(x,0)=u_{2}^{0}(x),\;u_{2}'(x,0)=u_{2}^{1}(x)&\text{on }(l,L),
\end{array}\right.
\end{equation}
and
\begin{equation}\label{kv2}
\left\{\begin{array}{ll}
\ddot{w}_{1}(x,t)+(q_{1}w_{1}''+D_{b}\dot{w}_{1}'')''(x,t)=0&\text{in }(0,l)\times(0,+\infty)
\\
\ddot{w}_{2}(x,t)+(q_{2}w_{2}'')''(x,t)=0&\text{in }(l,L)\times(0,+\infty)
\\
w_{1}(l,t)=w_{2}(l,t)&\text{for }t\in(0,+\infty)
\\
w_{1}'(l,t)=w_{2}'(l,t)&\text{for }t\in(0,+\infty)
\\
q_{1}.w_{1}''(l,t)=q_{2}.w_{2}''(l,t)&\text{for }t\in(0,+\infty)
\\
q_{1}.w_{1}'''(l,t)=q_{2}.w_{2}'''(l,t)&\text{for }t\in(0,+\infty)
\\
w_{1}(0,t)=w_{1}'(0,t)=w_{2}(L,t)=w_{2}'(L,t)=0&\text{for }t\in(0,+\infty)
\\
w_{1}(x,0)=w_{1}^{0}(x),\;w_{1}'(x,0)=w_{1}^{1}(x)&\text{on }(0,l)
\\
w_{2}(x,0)=w_{2}^{0}(x),\;w_{2}'(x,0)=w_{2}^{1}(x)&\text{on }(l,L),
\end{array}\right.
\end{equation}
where $u_{1}.\chi_{(0,l)}+u_{2}.\chi_{(l,L)}$ and $w_{1}.\chi_{(0,l)}+w_{2}.\chi_{(l,L)}$ represent the longitudinal and transversal displacement of the beam in the interval $(0,L)$ respectively,  with $\chi_{(l_{1},l_{2})}$ being the characteristic function of the interval $(l_{1},l_{2})$. The coefficient functions $p_{1},\,q_{1},\,D_{a},\,D_{b}$ are in $L^{\infty}(0,l)$, $p_{2},\,q_{2}$ are in $L^{\infty}(l,L)$ such that $p_{1},\,p_{2},\,q_{1},\,q_{2}\geq c_{0}>0$ and $D_{a}=a(x)\chi_{(\alpha,\beta)}$, $D_{b}=b(x)\chi_{(\alpha,\beta)}$, $0\leq \alpha<\beta\leq l<L$ and $a(x),\,b(x)\geq c_{0}>0$.

The energy of a solution of~\eqref{kv1} and~\eqref{kv2} at the time $t\geq 0$ are defined  respectively by
\[
E_{1}(t)=\frac{1}{2}\left(\int_{0}^{l}\Big(|\dot{u}_{1}(x,t)|^{2}+p_{1}|u_{1}'(x,t)|^{2}\Big)\,\ud x+\int_{l}^{L}\Big(|\dot{u}_{2}(x,t)|^{2}+p_{2}|u_{2}'(x,t)|^{2}\Big)\,\ud x\right).
\]
and
\[
E_{2}(t)=\frac{1}{2}\left(\int_{0}^{l}\Big(|\dot{w}_{1}(x,t)|^{2}+q_{1}|w_{1}''(x,t)|^{2}\Big)\,\ud x+\int_{l}^{L}\Big(|\dot{w}_{2}(x,t)|^{2}+q_{2}|w_{2}''(x,t)|^{2}\Big)\,\ud x\right).
\]
By Green's formula we can prove that for all $t_{1},\,t_{2}>0$ we have
\[
E_{1}(t_{2})-E_{1}(t_{1})=-\int_{t_{1}}^{t_{2}}\!\!\!\int_{\alpha}^{\beta}D_{a}|\dot{u}_{1}'(x,t)|^{2}\,\ud x\,\ud t,
\]
and
\[
E_{2}(t_{2})-E_{2}(t_{1})=-\int_{t_{1}}^{t_{2}}\!\!\!\int_{\alpha}^{\beta}D_{b}|\dot{w}_{1}''(x,t)|^{2}\,\ud x\,\ud t,
\]
this mean that the energy is decreasing over the time. In~\cite{LL} longitudinal and transversal vibrations of a clamped elastic beam were studied as problems with locally distributed damping. It was shown that when viscoelastic damping is distributed only on a subinterval in the interior of the domain, the exponential stability holds for the transversal but not for the longitudinal motion. Then it was proved in~\cite{LR1} that we have exactly a polynomial stability for the longitudinal motion.

At this point our main concern is the following question: Is the locally distributed Kelvin-Voigt damping (on any subinterval of $(0,l)$) strong enough to cause uniform exponential decay of the energy of the beam for transversal and longitudinal motion in the case of transmission problem given respectively by~\eqref{kv2} and~\eqref{kv1}?

An exponential stability of transmission problem for waves with frictional damping was treated in~\cite{BR}, for the Timoshenko system it was treated in~\cite{R1}, the general decay of solution for the transmission problem of viscoelastic waves with memory was treated in~\cite{R2} and in~\cite{LR2} uniform stability is proved for the wave equation with smooth viscoelastic damping applied just around the boundary and more recently in~\cite{BRA} in which the viscosity is distributed uniformly in the whole beam $(0,l)$ and in~\cite{ARSV} in which we consider a model of a material composed by three components, one of them is a Kelvin-Voigt viscoelastic material, the second is an elastic material (no dissipation) and the third is an elastic material inserted with a frictional damping mechanism.

Our main tool to prove the exponential stability 
is a result due to Pr\"uss~\cite{Pr} and to show the polynomial decay and the optimality of the decay rate we use a result due to Borichev and Tomilov~\cite{BT}.

The remaining part of this paper is organized as follows. In section~\ref{kv94} under some assumptions in the coefficients we prove that the energy decay of problem~\eqref{kv2} is exponentially stable and in section~\ref{kv95} we interested to the longitudinal motion given by~\eqref{kv1} for which we prove that the corresponding semigroup is polynomially stable and not exponentially stable.
\section{Transversal motion}\label{kv94}
Let $H=L^{2}(0,l)\times L^{2}(l,L)$ with the norm
$$\|v\|_{H}=\|(v_{1},v_{2})\|_{H}=\left(\int_{0}^{l}|v_{1}|^{2}(x)\,\ud x+\int_{l}^{L}|v_{2}|^{2}(x)\,\ud x\right)^{\frac{1}{2}}$$
and 
\begin{equation*}
\begin{split}
V_{1}=\{w=(w_{1},w_{2})\in H^{2}(0,l)\times H^{2}(l,L): w_{1}(l)=w_{2}(l),\,w_{1}'(l)=w_{2}'(l),
\\
w_{1}(0)=0,\,w_{2}(L)=0,\,w_{1}'(0)=0,\,w_{2}'(L)=0\}
\end{split}
\end{equation*}
with the norm
$$\|w\|_{V_{1}}=\|(w_{1},w_{2})\|_{V_{1}}=\left(\int_{0}^{l}q_{1}.|w_{1}''|^{2}(x)\,\ud x+\int_{l}^{L}q_{2}.|w_{2}''|^{2}(x)\,\ud x\right)^{\frac{1}{2}}.$$
Define $\mathcal{H}_{1}=V_{1}\times H$ with the norm $\|(w,v)\|_{\mathcal{H}_{1}}^{2}=\|w\|_{V_{1}}^{2}+\|v\|_{H}^{2}$. Then $\mathcal{H}_{1}$ is a Hilbert space in which we define
\begin{equation*}
\begin{split}
\mathcal{D}(\mathcal{A}_{1})=\{(w,v)\in\mathcal{H}_{1}:v\in V_{1},\,q_{1}.w_{1}''+D_{b}.v_{1}''\in H^{2}(0,l),\,w_{2}''\in H^{2}(l,L),
\\
q_{1}(l).w_{1}''(l)=q_{2}(l).w_{2}''(l),\,q_{1}(l).w_{1}'''(l)=q_{2}(l).w_{2}'''(l)\}
\end{split}
\end{equation*}
and
$$\mathcal{A}_{1}(w,v)=\mathcal{A}_{1}((w_{1},w_{2}),(v_{1},v_{2}))=(v_{1},v_{2},-(q_{1}.w_{1}''+D_{b}.v_{1}'')'',(q_{2}.w_{2}'')'').$$
Thus,~\eqref{kv2} can rewritten as an abstract evolution equation on $\mathcal{H}_{1}$,
$$(\dot{w}(t),\dot{v}(t))=\mathcal{A}_{1}(w(t),v(t)),\quad (w_{1}(0),w_{2}(0),v_{1}(0),v_{2}(0))=(w_{1}^{0},w_{2}^{0},w_{1}^{1},w_{2}^{1}).$$
\begin{pro}\label{kv96}
The linear operator $\mathcal{A}_{1}$ generates a $C_{0}$-semigroup of contractions $e^{t\mathcal{A}_{1}}$ on $\mathcal{H}_{1}$, in particular there exists a unique solution of~\eqref{kv2} which can be expressed by means of a semigroup on $\mathcal{H}_{1}$ having the following regularity of the solution
$$
\left(\begin{array}{l}
w_{1}
\\
w_{2}
\\
\dot{w}_{1}
\\
\dot{w}_{2}
\end{array}\right)
\in C([0,+\infty[,\mathcal{D}(\mathcal{A}_{1}))\cap C^{1}([0,+\infty[,\mathcal{H}_{1}),$$
if $(w_{1}^{0},w_{2}^{0},w_{1}^{1},w_{2}^{1})\in \mathcal{D}(\mathcal{A}_{1})$ and a mild solution
$$
\left(\begin{array}{l}
w_{1}
\\
w_{2}
\\
\dot{w}_{1}
\\
\dot{w}_{2}
\end{array}\right)
\in C([0,+\infty[,\mathcal{H}_{1}),$$
if $(w_{1}^{0},w_{2}^{0},w_{1}^{1},w_{2}^{1})\in\mathcal{H}_{1}$. In addition the imaginary axis $i\R$ is a subset of the resolvent set $\rho(\mathcal{A}_{1})$.
\end{pro}
\begin{pr}
The well-posedness problem follow easily from~\cite[section 2]{CLL} in which we have only to verify a simple conditions named (H1), (H2) and (H3), to prove the dissipative character of the operator, which is the case for our problem.
\\
It is easy to show that there is no point spectrum of $\mathcal{A}_{1}$ on the imaginary axis, i.e. $i\R\cap\sigma_{p}(\mathcal{A}_{1})\neq\emptyset$. Further, $V_{1}\hookrightarrow H$ with a compact embedding then the result of the resolvent set follow from~\cite[Lemma 4.1]{CLL}.
\end{pr}

We assume that $q_{1}$ is constant in each of the interval $(0,\alpha)$, $(\beta,l)$ and $q_{2}$ is also constant in $(l,L)$ and $b,\,q_{1}\in \mathrm{C}([\alpha,\beta])$, and we suppose that we have
\begin{equation}\label{kv42}
q_{2}(l)\geq q_{1}(l).
\end{equation}
\begin{thm}
Under the above assumptions on the coefficients of~\eqref{kv2}, the semigroup $e^{\mathcal{A}_{1}t}$ is exponentially stable, i.e., there exist $\mu>0$ and $C_{\mu}>0$ such that
$$\|e^{\mathcal{A}_{1}t}\|_{\mathcal{H}}\leq C_{\mu}\e^{-\mu t},\quad\forall\,t\geq0.$$
\end{thm}
\begin{pr}
We need only to verify the condition for a $C_{0}$-semigroup of contractions on a Hilbert space being exponentially stable (see~\cite{Hu},~\cite{Pr} or~\cite{G}) i.e.,
\begin{equation}\label{kv4}
\sup\{\|(i\lambda-\mathcal{A}_{1})^{-1}\|:\lambda\in\R\}<+\infty.
\end{equation}
Suppose that~\eqref{kv4} is not true. By the continuity of the resolvent and the resonance theorem, there exist $\lambda_{n}\in\R$, $(w_{n},v_{n})=((w_{1,n},w_{2,n}),(v_{1,n},v_{2,n}))\in\mathcal{D}(\mathcal{A}_{1})$, for all $n\in\N$ such that
\begin{equation}\label{kv3}
\|(w_{n},v_{n})\|_{\mathcal{H}_{1}}=1,\;\lim_{n\rightarrow+\infty}|\lambda_{n}|=+\infty
\end{equation}
and
\begin{equation}\label{kv15}
(i\lambda_{n}-\mathcal{A}_{1})(w_{n},v_{n})=(f_{n},g_{n})=((f_{1,n},f_{2,n}),(g_{1,n},g_{2,n}))\,\longrightarrow\, 0\text{ in }\mathcal{H}_{1}.
\end{equation}
This implies
\begin{eqnarray}
i\lambda_{n}w_{1,n}-v_{1,n}=f_{1,n}\,\longrightarrow\,0\text{ in }H^{2}(0,l),\label{kv5}
\\
i\lambda_{n}w_{2,n}-v_{2,n}=f_{2,n}\,\longrightarrow\,0\text{ in }H^{2}(l,L),\label{kv6}
\\
i\lambda_{n}.v_{1,n}-M_{n}''=g_{1,n}\,\longrightarrow\,0\text{ in }L^{2}(0,l),\label{kv7}
\\
i\lambda_{n}.v_{2,n}+q_{2}w_{2,n}''''=g_{2,n}\,\longrightarrow\,0\text{ in }L^{2}(l,L),\label{kv8}
\end{eqnarray}
where 
$$
M_{n}=-(q_{1}w_{1,n}''+D_{b}v_{1,n}'').
$$
For $0\leq x\leq l$ define
\begin{equation}\label{kv9}
J(\psi)(x)=\int_{x}^{l}\!\!\int_{s}^{l}\psi(\tau)\,\ud\tau\,\ud s
\end{equation}
and
\begin{equation}\label{kv10}
y_{n}=\frac{1}{i\lambda_{n}}[M_{n}+J(g_{1,n})].
\end{equation}
Comparing~\eqref{kv9} and~\eqref{kv10} we have
\begin{equation}\label{kv11}
y_{n}''=v_{1,n}.
\end{equation}
The rest of the proof depends on the following two lemmas. Let $\omega_{n}=\sqrt{|\lambda_{n}|}$.
\begin{lem}
The function $y_{n}$ defined above has the following properties:
\begin{align}
y_{n}\,\longrightarrow\,0\quad\text{ in }\; H^{4}(\alpha,\beta),\label{kv12}
\\
\lambda_{n}.y_{n}\,\longrightarrow\,0\quad\text{ in }\; L^{2}(\alpha,\beta),\label{kv13}
\\
\omega_{n}.y_{n}\,\longrightarrow\,0\quad\text{ in }\; H^{2}(\alpha,\beta).\label{kv14}
\end{align}
\end{lem}
\begin{pr}
From~\eqref{kv15},
\begin{equation}\label{kv16}
\re\langle(i\lambda_{n}-\mathcal{A}_{1})(w_{n},v_{n}),(w_{n},v_{n})\rangle_{\mathcal{H}_{1}}=\int_{\alpha}^{\beta}D_{b}|v_{1,n}''|^{2}\,\ud x\,\longrightarrow\,0.
\end{equation}
Therefore, from~\eqref{kv5} we have
\begin{equation}\label{kv17}
M_{n}\,\longrightarrow\,0\quad\text{ in }\;L^{2}(\alpha,\beta)
\end{equation}
and
\begin{equation}\label{kv18}
\frac{1}{|\lambda_{n}|}\|\psi.v_{n}\|_{V_{1}}=O(1)
\end{equation}
for every $\psi\in\Ci([0,l])$, such that $\supp(\psi)\subset[0,l]$.

Equations~\eqref{kv7},~\eqref{kv17} and~\eqref{kv18} imply that
\begin{equation}\label{kv19}
\int_{\alpha}^{\beta}\psi.|v_{1,n}|^{2}\,\ud x\,\longrightarrow\,0,\quad\forall\,\psi\in\Ci([0,l]),\;\supp(\psi)\subset(\alpha,\beta).
\end{equation}
Applying the interpolation theorem involving compact subdomain~\cite[Theorem 4.23]{A} we find that~\eqref{kv16} and~\eqref{kv19} imply
\begin{equation}\label{kv20}
v_{1,n}\,\longrightarrow\,0\quad\text{ in }\;H^{2}(\alpha,\beta)
\end{equation}
Thus,~\eqref{kv11} yields
\begin{equation}\label{kv21}
\int_{\alpha}^{\beta}|y_{n}''''|^{2}\,\ud x\,\longrightarrow\,0.
\end{equation}
In the other hand,~\eqref{kv13} follow from
\begin{equation}\label{kv22}
i\lambda_{n}y_{n}=J(g_{1,n})-\frac{q_{1}+i\lambda_{n}D_{b}}{i\lambda_{n}}v_{n}''-\frac{q_{1}}{i\lambda_{n}}f_{1,n}''\,\longrightarrow\,0\quad\text{ in }\;L^{2}(\alpha,\beta).
\end{equation}
Since $|\lambda_{n}|\,\longrightarrow\,+\infty$, we obtain that $y_{n}\,\longrightarrow\,0$ in $L^{2}(\alpha,\beta)$. This combined with~\eqref{kv21} yields~\eqref{kv12}. From the interpolation inequality~\cite[Theorem 4.17]{A}, we also have~\eqref{kv14}.
\end{pr}
\begin{lem}\label{kv111}
The functions $w_{1,n}\in H^{4}(0,\alpha)\cap H^{4}(\beta,l)$, for all $n\in\N$ have the following properties:
\begin{align}
\omega_{n}^{4}(|w_{1,n}(\alpha)|^{2}+|w_{1,n}'(\alpha)|^{2}+|w_{1,n}(\beta)|^{2}+|w_{1,n}'(\beta)|^{2})\,\longrightarrow\,0,\label{kv23}
\\
\alpha q_{1}(0)|w_{1,n}''(\alpha^{-})|^{2}+(L-\beta)q_{1}(l)|w_{1,n}''(\beta^{+})|^{2}+(L-l)q_{1}(l)\left(\frac{q_{1}(l)}{q_{2}(l)}-1\right)|w_{1,n}''(l)|^{2}\longrightarrow\,2,\label{kv24}
\\
\omega_{n}^{-1}|w_{1,n}'''(\alpha^{-})|\,\longrightarrow\,0\quad\text{ and }\quad\omega_{n}^{-1}|w_{1,n}'''(\beta^{+})|\,\longrightarrow\,0.\label{kv25}
\end{align}
\end{lem}
\begin{pr}
Since $w_{1,n},\,v_{1,n}\in H^{2}(0,l)$, Sobolev's embedding theorem implies that they are also in $C^{1}(0,l)$. By~\eqref{kv5} and~\eqref{kv20} we have
\begin{equation}\label{kv26}
\lambda_{n}.w_{1,n}\,\longrightarrow\,0\quad\text{ in }\;H^{2}(\alpha,\beta).
\end{equation}
Thus, $\lambda_{n}.w_{1,n}$ converges to zero in $C^{1}([\alpha,\beta])$, which immediately leads to~\eqref{kv23}.

Note that $M_{n}=-q_{1}w_{1,n}''$ on $(0,\alpha)\cup(\beta,l)$, where $q_{1}\equiv q_{1}(0)$ on $[0,\alpha)$ and $q_{1}\equiv q_{1}(l)$ on $(\beta,l]$. From the definition of the domain of $\mathcal{A}_{1}$, we know $w_{1,n}\in H^{4}(0,\alpha)$ and $w_{1,n}\in H^{4}(\beta,l)$. It follows from~\eqref{kv10} that
\begin{align}
q_{1}(0).w_{1,n}''(\alpha^{-})=(J(g_{1,n})-i\lambda_{n}y_{n})(\alpha),\qquad q_{1}(l).w_{1,n}''(\beta^{+})=(J(g_{1,n})-i\lambda_{n}y_{n})(\beta),\label{kv27}
\\
q_{1}(0).w_{1,n}'''(\alpha^{-})=(J(g_{1,n})-i\lambda_{n}y_{n})'(\alpha),\qquad q_{1}(l).w_{1,n}'''(\beta^{+})=(J(g_{1,n})-i\lambda_{n}y_{n})'(\beta).\label{kv28}
\end{align}
Dividing~\eqref{kv28} by $\omega_{n}$ we obtain~\eqref{kv25} by using~\eqref{kv14} in the previous lemma.

In order to prove~\eqref{kv24}, we substitute~\eqref{kv5} into~\eqref{kv7} and~\eqref{kv6} into~\eqref{kv8} to get
\begin{equation}\label{kv29}
\left\{\begin{array}{ll}
-\lambda_{n}^{2}.w_{1,n}-M_{n}''=g_{1,n}+i\lambda_{n}.f_{1,n}&\text{in }(0,l)
\\
-\lambda_{n}^{2}.w_{2,n}+q_{2}.w_{2,n}''''=g_{2,n}+i\lambda_{n}.f_{2,n}&\text{in }(l,L).
\end{array}\right.
\end{equation}
We multiply the above equations by $\overline{w}_{1,n}$ and $\overline{w}_{2,n}$ respectively, then integrate by parts on $(0,L)$. This leads to
\begin{equation}\label{kv30}
\|\lambda_{n}.w_{n}\|_{H}^{2}-\|w_{n}\|_{V_{1}}^{2}\,\longrightarrow\,0.
\end{equation}
Here we have used~\eqref{kv3},~\eqref{kv5},~\eqref{kv6},~\eqref{kv7},~\eqref{kv8} and~\eqref{kv17}. Since $\|w_{n}\|_{V_{1}}^{2}+\|v_{n}\|_{H}^{2}=1$ and $i\lambda_{n}.w_{n}-v_{n}$ also converges to zero in $L^{2}(0,L)$,~\eqref{kv30} implies that both $\|\lambda_{n}.w_{n}\|_{H}^{2}$ and $\|w_{n}\|_{V_{1}}^{2}$ must converge to $\displaystyle\frac{1}{2}$ as $n\,\longrightarrow\,+\infty$. This further leads to
\begin{equation}\label{kv31}
\begin{split}
\lim_{n\rightarrow+\infty}\int_{0}^{\alpha}|\lambda_{n}.w_{1,n}|^{2}\,\ud x+\int_{\beta}^{l}|\lambda_{n}.w_{1,n}|^{2}\,\ud x+\int_{l}^{L}|\lambda_{n}.w_{2,n}|^{2}\,\ud x=
\\
\lim_{n\rightarrow+\infty}\int_{0}^{\alpha}q_{1}|w_{1,n}''|^{2}\,\ud x+\int_{\beta}^{l}q_{1}|w_{1,n}''|^{2}\,\ud x+\int_{l}^{L}q_{2}|w_{2,n}''|^{2}\,\ud x=\frac{1}{2}
\end{split}
\end{equation}
when~\eqref{kv26} is taken into account.

On the intervals $(0,\alpha)$, $(\beta,l)$ and $(l,L)$,~\eqref{kv29} becomes
\begin{equation}\label{kv32}
\left\{\begin{array}{ll}
-\lambda_{n}^{2}.w_{1,n}+q_{1}.w_{1,n}''''=g_{1,n}+i\lambda_{n}.f_{1,n}&\text{in }(0,\alpha)\cup(\beta,l)
\\
-\lambda_{n}^{2}.w_{2,n}+q_{2}.w_{2,n}''''=g_{2,n}+i\lambda_{n}.f_{2,n}&\text{in }(l,L).
\end{array}\right.
\end{equation}
We multiply the first equation of~\eqref{kv32} respectively by $x\overline{w'}_{1,n}$ and $(L-x)\overline{w'}_{1,n}$ and integrate on $(0,\alpha)$ and $(\beta,l)$ respectively. Hence,
\begin{equation}\label{kv33}
-\int_{0}^{\alpha}\lambda_{n}^{2}w_{1,n}x\overline{w'}_{1,n}\,\ud x+\int_{0}^{\alpha}q_{1}w_{1,n}''''x\overline{w'}_{1,n}\,\ud x=\int_{0}^{\alpha}(g_{1,n}+i\lambda_{n}f_{1,n})x\overline{w'}_{1,n}\,\ud x
\end{equation}
and
\begin{equation}\label{kv34}
\begin{split}
-\int_{\beta}^{l}\lambda_{n}^{2}w_{1,n}(L-x)\overline{w'}_{1,n}\,\ud x+\int_{\beta}^{l}q_{1}w_{1,n}''''(L-x)\overline{w'}_{1,n}\,\ud x
\\
=\int_{\beta}^{l}(g_{1,n}+i\lambda_{n}f_{1,n})(L-x)\overline{w'}_{1,n}\,\ud x.
\end{split}
\end{equation}
It is easy to see that the terms on the right hand side of~\eqref{kv33} and~\eqref{kv34} converges to zero. After a straightforward calculation (integration by parts), the two terms on the left hand side of~\eqref{kv33} and~\eqref{kv34} become
\begin{equation}\label{kv35}
-\re\left(\int_{0}^{\alpha}\lambda_{n}^{2}w_{1,n}x\overline{w'}_{1,n}\,\ud x\right)=-\omega_{n}^{4}\frac{\alpha}{2}|w_{1,n}(\alpha)|^{2}+\frac{1}{2}\int_{0}^{\alpha}|\lambda_{n}w_{1,n}|^{2}\,\ud x,
\end{equation}
\begin{equation}\label{kv36}
\begin{split}
\re\left(\int_{0}^{\alpha}q_{1}w_{1,n}''''x\overline{w'}_{1,n}\,\ud x\right)=q_{1}(0)\re[(\alpha w_{1,n}'''(\alpha^{-})-w_{1,n}''(\alpha^{-}))\overline{w'}_{1,n}(\alpha)]
\\
+\frac{3}{2}\int_{0}^{\alpha}q_{1}|w_{1,n}''|^{2}\,\ud x-\frac{\alpha}{2}q_{1}(0)|w_{1,n}''(\alpha^{-})|^{2},
\end{split}
\end{equation}
and
\begin{equation}\label{kv37}
\begin{split}
-\re\left(\int_{\beta}^{l}\lambda_{n}^{2}w_{1,n}(L-x)\overline{w'}_{1,n}\,\ud x\right)=-\frac{1}{2}\int_{\beta}^{l}|\lambda_{n}w_{1,n}|^{2}\,\ud x
\\
+\frac{\omega_{n}^{4}}{2}[(L-\beta)|w_{1,n}(\beta)|^{2}-(L-l)|w_{1,n}(l)|^{2}],
\end{split}
\end{equation}
\begin{equation}\label{kv38}
\begin{split}
\re\left(\int_{\beta}^{l}q_{1}w_{1,n}''''(L-x)\overline{w'}_{1,n}\,\ud x\right)=q_{1}(l)\re[((L-l) w_{1,n}'''(l)+w_{1,n}''(l))\overline{w'}_{1,n}(l)]
\\
-q_{1}(l)\re[((L-\beta) w_{1,n}'''(\beta^{+})+w_{1,n}''(\beta^{+}))\overline{w'}_{1,n}(\beta)]-\frac{3}{2}\int_{\beta}^{l}q_{1}|w_{1,n}''|^{2}\,\ud x
\\
+\frac{(L-\beta)}{2}q_{1}(l)|w_{1,n}''(\beta^{+})|^{2}-\frac{(L-l)}{2}q_{1}(l)|w_{1,n}''(l)|^{2}.
\end{split}
\end{equation}
After substituting these terms into the real part of~\eqref{kv33} and~\eqref{kv34} and applying~\eqref{kv23},~\eqref{kv25} and~\eqref{kv26}, we obtain
\begin{equation}\label{kv39}
\frac{1}{2}\int_{0}^{\alpha}|\lambda_{n}w_{1,n}|^{2}\,\ud x+\frac{3}{2}\int_{0}^{\alpha}q_{1}|w_{1,n}''|^{2}\,\ud x-\frac{\alpha}{2}q_{1}(0)|w_{1,n}''(\alpha^{-})|^{2}\,\longrightarrow\,0,
\end{equation}
and
\begin{equation}\label{kv40}
\begin{split}
\frac{1}{2}\int_{\beta}^{l}|\lambda_{n}w_{1,n}|^{2}\,\ud x+\frac{3}{2}\int_{\beta}^{l}q_{1}|w_{1,n}''|^{2}\,\ud x-q_{1}(l)\re[((L-l)w_{1,n}'''(l)+w_{1,n}''(l))\overline{w'}_{1,n}(l)]
\\
-\frac{(L-\beta)}{2}q_{1}(l)|w_{1,n}''(\beta^{+})|^{2}+\frac{(L-l)}{2}q_{1}(l)|w_{1,n}''(l)|^{2}+\frac{(L-l)}{2}\omega_{n}^{4}|w_{1,n}(l)|^{2}\,\longrightarrow\,0.
\end{split}
\end{equation}
Similarly, we can multiply the second equation of~\eqref{kv32} by $(L-x)\overline{w'}_{2,n}$ and integrate on $(l,L)$ to get
\begin{equation}\label{kv41}
\begin{split}
\frac{1}{2}\int_{l}^{L}|\lambda_{n}w_{2,n}|^{2}\,\ud x+\frac{3}{2}\int_{l}^{L}q_{2}|w_{2,n}''|^{2}\,\ud x+q_{2}(l)\re[((L-l)w_{2,n}'''(l)+w_{2,n}''(l))\overline{w'}_{2,n}(l)]
\\
-\frac{(L-l)}{2}q_{2}(l)|w_{2,n}''(l)|^{2}-\frac{(L-l)}{2}\omega_{n}^{4}|w_{2,n}(l)|^{2}\,\longrightarrow\,0.
\end{split}
\end{equation}
Then by summing~\eqref{kv39},~\eqref{kv40} and~\eqref{kv41} and using~\eqref{kv31} and the transmission conditions we find that
$$
\alpha q_{1}(0)|w_{1,n}''(\alpha^{-})|^{2}+(L-\beta)q_{1}(l)|w_{1,n}''(\beta^{+})|^{2}+(L-l)(q_{2}(l)|w_{2,n}''(l)|^{2}-q_{1}(l)|w_{1,n}''(l)|^{2})\longrightarrow\,2.
$$
Finally, by using again the transmission conditions we obtain~\eqref{kv24} easly.
\end{pr}

In what follows, we will show that
\begin{equation}\label{kv43}
|w_{1,n}''(\alpha^{-})|^{2}\,\longrightarrow\,0\quad\text{ and }\quad |w_{1,n}''(\beta^{+})|^{2}\,\longrightarrow\,0,
\end{equation}
to obtain a contradiction by the assumption~\eqref{kv42} with~\eqref{kv24}. Denote by
\begin{equation*}
\begin{array}{lcll}
\displaystyle\fo=\frac{\omega_{n}}{\sqrt[4]{q_{1}}},&\displaystyle\ft=\frac{\omega_{n}}{\sqrt[4]{q_{2}}},&\displaystyle\fot=\frac{\fo}{\ft},&\displaystyle\fto=\frac{\ft}{\fo},
\\
\displaystyle F_{1,n}=\frac{1}{q_{1}}(g_{1,n}+i\lambda_{n} f_{1,n})&\text{and}&\displaystyle F_{2,n}=\frac{1}{q_{2}}(g_{2,n}+i\lambda_{n} f_{2,n}).
\end{array}
\end{equation*}
Then~\eqref{kv32} can be writen as
\begin{equation}\label{kv44}
\left\{\begin{array}{ll}
\displaystyle(\de-i\fo)(\de+i\fo)(\de^{2}-\fo^{2})w_{1,n}=F_{1,n}&\text{ in }(0,\alpha)\cup(\beta,l)
\\
\displaystyle(\de-i\ft)(\de+i\ft)(\de^{2}-\ft^{2})w_{2,n}=F_{2,n}&\text{ in }(l,L),
\end{array}\right.
\end{equation}
where  we denoted by $\displaystyle\de=\frac{\ud}{\ud x}$.

The main idea consist to solve the equation~\eqref{kv44} in each of the intervals $(0,\alpha),\,(\beta,l)$ and $(l,L)$ in several steps, then using the boundary and transmission conditions we obtain informations about the traces of $w_{1,n}$ in $x=\alpha$ and $x=\beta$ which will not be compatible with those of Lemma~\ref{kv111}.

On the interval $(0,\alpha)$, by solving the first linear equation of~\eqref{kv44} we get
\begin{equation*}
(\de+i\fo)(\de^{2}-\fo^{2})w_{1,n}=C_{1}'\e^{i\fo(x-\alpha)}+\int_{\alpha}^{x}\e^{i\fo(x-s)}F_{1,n}(s)\,\ud s,
\end{equation*}
and
\begin{equation}\label{kv46}
(\de^{2}-\fo^{2})w_{1,n}=C_{2}'\e^{-i\fo(x-\alpha)}+\frac{C_{1}'}{\fo}\sin(\fo(x-\alpha))+\frac{1}{\fo}\int_{\alpha}^{x}\sin(\fo(x-s))F_{1,n}(s)\,\ud s,
\end{equation}
where
\begin{equation}\label{kv51}
C_{1}'=w_{1,n}'''(\alpha^{-})+i\fo w_{1,n}''(\alpha^{-})-\fo^{2}w_{1,n}'(\alpha)-i\fo^{3}w_{1,n}(\alpha),
\end{equation}
and
\begin{equation}\label{kv52}
C_{2}'=w_{1,n}''(\alpha^{-})-\fo^{2}w_{1,n}(\alpha).
\end{equation}
We still resolves again~\eqref{kv46} and using the boundary conditions $w_{1,n}(0)=w_{1,n}'(0)=0$ to get
\begin{equation}\label{kv48}
\begin{split}
(\de-\fo)w_{1,n}=\frac{C_{2}'}{(1-i)\fo}\left[\e^{-i\fo(x-\alpha)}-\e^{-\fo(x-i\alpha)}\right]
\\
+\frac{C_{1}'}{2\fo^{2}}\left[\sin(\fo(x-\alpha))-\cos(\fo(x-\alpha))+\e^{-\fo x}(\sin(\fo\alpha)+\cos(\fo\alpha))\right]
\\
+\frac{1}{\fo}\int_{0}^{x}\int_{\alpha}^{\tau}\e^{-\fo(x-\tau)}\sin(\fo(\tau-s))F_{1,n}(s)\,\ud s\,\ud\tau.
\end{split}
\end{equation}
Multiplying~\eqref{kv48} by $2\fo$ and taking $x=\alpha$, we have
\begin{equation}\label{kv53}
\begin{split}
2\fo w_{1,n}'(\alpha)-2\fo^{2}w_{1,n}(\alpha)=(1+i)C_{2}'(1-\e^{-\alpha\fo}\e^{i\alpha\fo})
\\
+\frac{C_{1}'}{\fo}(\e^{-\alpha\fo}(\sin(\alpha\fo)+\cos(\alpha\fo))-1)
\\
+2\int_{0}^{\alpha}\int_{\alpha}^{\tau}\e^{-\fo(\alpha-\tau)}\sin(\fo(\tau-s))F_{1,n}(s)\,\ud s\,\ud\tau.
\end{split}
\end{equation}
We substitute~\eqref{kv51} and~\eqref{kv52} into~\eqref{kv53} and let $n\,\longrightarrow\,+\infty$, then~\eqref{kv23} and~\eqref{kv25} yields
$$\lim_{n\rightarrow+\infty}w_{1,n}''(\alpha^{-})=-2\lim_{n\rightarrow+\infty}\int_{0}^{\alpha}\int_{\alpha}^{\tau}\e^{-\fo(\alpha-\tau)}\sin(\fo(\tau-s))F_{1,n}(s)\,\ud s\,\ud\tau.$$
We argue that the above limit is zero by the following estimates
$$\left|\int_{0}^{\alpha}\int_{\alpha}^{\tau}\e^{-\fo(\alpha-\tau)}\sin(\fo(\tau-s))g_{1,n}(s)\,\ud s\,\ud\tau\right|\leq\alpha^{\frac{3}{2}}\left(\int_{0}^{\alpha}|g_{1,n}(s)|^{2}\ud s\right)^{\frac{1}{2}}\,\longrightarrow\,0,$$
and
\begin{equation*}
\begin{split}
\left|\int_{0}^{\alpha}\int_{\alpha}^{\tau}\e^{-\fo(\alpha-\tau)}\sin(\fo(\tau-s))\lambda_{n}f_{1,n}(s)\,\ud s\,\ud\tau\right|=
\\
\left|\lambda_{n}\int_{0}^{\alpha}\left(\int_{0}^{s}\e^{-\fo(\alpha-\tau)}\sin(\fo(\tau-s))\,\ud\tau\right)f_{1,n}(s)\,\ud s\right|=
\\
\left|\frac{\lambda_{n}\e^{-\alpha\fo}}{2\fo}\int_{0}^{\alpha}(\cos(\fo s)+\sin(\fo s)-\e^{\fo s})f_{1,n}(s)\,\ud s\right|\leq
\\
\sqrt{q_{1}}\left(\alpha\fo\e^{-\alpha\fo}+\frac{1}{2}-\frac{1}{2}\e^{-\alpha\fo}\right)\max_{[0,\alpha]}|f_{1,n}(s)|\,\longrightarrow\,0,
\end{split}
\end{equation*}
where we have used the fact that $g_{n}\,\longrightarrow\,0$ in $H$, $f_{n}\,\longrightarrow\,0$ in $V_{1}\hookrightarrow\mathrm{C}^{1}([0,l])$, and $\fo\,\longrightarrow\,+\infty$. Thus we have proved the first identity of~\eqref{kv43}.

On the interval $(l,L)$, by solving the first linear equation~\eqref{kv44} we get
\begin{equation*}
(\de+i\ft)(\de^{2}-\ft^{2})w_{2,n}=K_{1}\e^{i\ft(x-l)}+\int_{l}^{x}\e^{i\ft(x-s)}F_{2,n}(s)\,\ud s,
\end{equation*}
\begin{equation*}
(\de^{2}-\ft^{2})w_{2,n}=K_{2}\e^{-i\ft(x-l)}+\frac{K_{1}}{\ft}\sin(\ft(x-l))+\frac{1}{\ft}\int_{l}^{x}\sin(\ft(x-s))F_{2,n}(s)\,\ud s,
\end{equation*}
where
\begin{eqnarray}\label{kv55}
K_{1}=w_{2,n}'''(l)+i\ft w_{2,n}''(l)-\ft^{2}w_{2,n}'(l)-i\ft^{3}w_{2,n}(l),
\end{eqnarray}
\begin{equation}\label{kv56}
K_{2}=w_{2,n}''(l)-\ft^{2}w_{2,n}(l)
\end{equation}
and using the boundary conditions $w_{2,n}(L)=w_{2,n}'(L)=0$ as previously then we get
\begin{equation}\label{kv50}
\begin{split}
(\de+\ft)w_{2,n}=\frac{K_{2}}{(1+i)\ft}\left[\e^{-\ft(L-x+i(L-l))}-\e^{-i\ft(x-l)}\right]
\\
+\frac{K_{1}}{2\ft^{2}}\left[-\cos(\ft(x-l))-\sin(\ft(x-l))+\e^{\ft (x-L)}(\cos(\ft(L-l))+\sin(\ft(L-l)))\right]
\\
+\frac{1}{\ft}\int_{L}^{x}\int_{l}^{\tau}\e^{\ft(x-\tau)}\sin(\ft(\tau-s))F_{2,n}(s)\,\ud s\,\ud\tau.
\end{split}
\end{equation}

On the interval $(\beta,l)$ similary as previously we have
\begin{equation}\label{kv45}
(\de+i\fo)(\de^{2}-\fo^{2})w_{1,n}=C_{1}\e^{i\fo(x-\beta)}+\int_{\beta}^{x}\e^{i\fo(x-s)}F_{1,n}(s)\,\ud s,
\end{equation}
\begin{equation}\label{kv47}
(\de^{2}-\fo^{2})w_{1,n}=C_{2}\e^{-i\fo(x-\beta)}+\frac{C_{1}}{\fo}\sin(\fo(x-\beta))+\frac{1}{\fo}\int_{\beta}^{x}\sin(\fo(x-s))F_{1,n}(s)\,\ud s,
\end{equation}
where
\begin{equation}\label{kv54}
C_{1}=w_{1,n}'''(\beta^{+})+i\fo w_{1,n}''(\beta^{+})-\fo^{2}w_{1,n}'(\beta)-i\fo^{3}w_{1,n}(\beta),
\end{equation}
\begin{equation}\label{kv92}
C_{2}=w_{1,n}''(\beta^{+})-\fo^{2}w_{1,n}(\beta),
\end{equation}
and
\begin{equation}\label{kv49}
\begin{split}
(\de+\fo)w_{1,n}=\frac{C_{2}}{(1+i)\fo}\left[\e^{\fo(x-\beta)}-\e^{-i\fo(x-\beta)}\right]
\\
+\frac{C_{1}}{2\fo^{2}}\left[1-\cos(\fo(x-\beta))-\sin(\fo(x-\beta))\right]+C_{3}\e^{\fo (x-\beta)}
\\
+\frac{1}{\fo}\int_{\beta}^{x}\int_{\beta}^{\tau}\e^{\fo(x-\tau)}\sin(\fo(\tau-s))F_{1,n}(s)\,\ud s\,\ud\tau,
\end{split}
\end{equation}
where
\begin{equation}\label{kv57}
C_{3}=w_{1,n}'(\beta)+\fo w_{1,n}(\beta).
\end{equation}
We resolve the first order equation of~\eqref{kv49} then we obtain
\begin{equation}\label{kv59}
\begin{split}
w_{1,n}=\frac{C_{2}}{2(1+i)\fo^{2}}\left[\e^{\fo(x-\beta)}+i\e^{-\fo(x-\beta)}-(1+i)\e^{-i\fo(x-\beta)}\right]
\\
+\frac{C_{1}}{2\fo^{3}}\left[1-\sin(\fo(x-\beta))-\e^{-\fo(x-\beta)}\right]
\\
+\frac{C_{3}}{2\fo}\left[\e^{\fo(x-\beta)}-\e^{-\fo(x-\beta)}\right]+C_{4}\e^{-\fo(x-\beta)}
\\
+\frac{1}{\fo}\int_{\beta}^{x}\int_{\beta}^{t}\int_{\beta}^{\tau}\e^{\fo(2t-\tau-x)}\sin(\fo(\tau-s))F_{1,n}(s)\,\ud s\,\ud\tau\,\ud t
\end{split}
\end{equation}
where
\begin{equation}\label{kv93}
C_{4}=w_{1,n}(\beta).
\end{equation}
Multiplying the relation~\eqref{kv59} by $\e^{-\fo(x-\beta)}$, taking $x=l$ and using the transmission condition $w_{1,n}(l)=w_{2,n}(l)$ then we get
\begin{equation}\label{kv60}
\begin{split}
\e^{-\fo(l-\beta)}w_{2,n}(l)=\e^{-\fo(l-\beta)}w_{1,n}(l)=
\\
\frac{C_{2}(1-i)}{4\fo^{2}}\left[1+i\e^{-2\fo(l-\beta)}-(1+i)\e^{-\fo(1+i)(l-\beta)}\right]
\\
+\frac{C_{1}}{2\fo^{3}}\left[\e^{-\fo(l-\beta)}-\sin(\fo(l-\beta))\e^{-\fo(l-\beta)}-\e^{-2\fo(l-\beta)}\right]
\\
+\frac{C_{3}}{2\fo}\left[1-\e^{-2\fo(l-\beta)}\right]+C_{4}\e^{-2\fo(l-\beta)}
\\
+\frac{1}{\fo}\int_{\beta}^{l}\int_{\beta}^{t}\int_{\beta}^{\tau}\e^{\fo(2t-\tau-2l+\beta)}\sin(\fo(\tau-s))F_{1,n}(s)\,\ud s\,\ud\tau\,\ud t,
\end{split}
\end{equation}
with integrating by parts the last term can be written as follows
\begin{equation}\label{kv107}
\begin{split}
&\int_{\beta}^{l}\int_{\beta}^{t}\int_{\beta}^{\tau}\e^{\fo(2t-\tau-2l+\beta)}\sin(\fo(\tau-s))F_{1,n}(s)\,\ud s\,\ud\tau\,\ud t=
\\
&\frac{1}{2\fo}\bigg(\int_{\beta}^{l}\e^{\fo(\beta-\tau)}\int_{\beta}^{\tau}\sin(\fo(\tau-s))F_{1,n}(s)\,\ud s\,\ud\tau
\\
&-\int_{\beta}^{l}\e^{\fo(t-2l+\beta)}\int_{\beta}^{t}\sin(\fo(t-s))F_{1,n}(s)\,\ud s\,\ud t\bigg)=
\\
&\frac{1}{2\fo}\bigg(\int_{\beta}^{l}\int_{s}^{l}\e^{\fo(\beta-\tau)}\sin(\fo(\tau-s))\,\ud \tau F_{1,n}(s)\,\ud s
\\
&-\int_{\beta}^{l}\int_{s}^{l}\e^{\fo(t-2l+\beta)}\sin(\fo(t-s))\,\ud t F_{1,n}(s)\,\ud s\bigg)=
\\
&\frac{1}{4\fo^{2}}\bigg(\int_{\beta}^{l}\e^{-\fo(s-\beta)}F_{1,n}(s)\,\ud s+\int_{\beta}^{l}\e^{-\fo(2l-\beta-s)}F_{1,n}(s)\,\ud s
\\
&-2\,\e^{-\fo(l-\beta)}\int_{\beta}^{l}\cos(\fo(l-s))F_{1,n}(s)\,\ud s\bigg),
\end{split}
\end{equation}
where the right hand side terms of~\eqref{kv107} verify
\begin{eqnarray}
\frac{1}{\fo}\int_{\beta}^{l}\e^{-\fo(s-\beta)}F_{1,n}(s)\,\ud s\,\longrightarrow\,0,\label{kv108}
\\
\frac{1}{\fo}\int_{\beta}^{l}\e^{-\fo(2l-\beta-s)}F_{1,n}(s)\,\ud s\,\longrightarrow\,0,\label{kv109}
\\
\frac{1}{\fo}\e^{-\fo(l-\beta)}\int_{\beta}^{l}\cos(\fo(l-s))F_{1,n}(s)\,\ud s\,\longrightarrow\,0,\label{kv110}
\end{eqnarray}
where we argue this by the fact that
\begin{equation*}
\begin{split}
\frac{1}{\fo}\left|\int_{\beta}^{l}\e^{-\fo(s-\beta)}g_{1,n}(s)\,\ud s\right|\leq\frac{(l-\beta)}{\fo}^{\frac{1}{2}}\left(\int_{\beta}^{l}|g_{1,n}(s)|^{2}\ud s\right)^{\frac{1}{2}}\,\longrightarrow\,0,
\\
\left|\int_{\beta}^{l}\e^{-\fo(2l-\beta-s)}g_{1,n}(s)\,\ud s\right|\leq\e^{-\fo(l-\beta)}(l-\beta)^{\frac{1}{2}}\left(\int_{\beta}^{l}|g_{1,n}(s)|^{2}\ud s\right)^{\frac{1}{2}}\,\longrightarrow\,0,
\\
\e^{-\fo(l-\beta)}\left|\int_{\beta}^{l}\cos(\fo(l-s))g_{1,n}(s)\,\ud s\right|\leq\e^{-\fo(l-\beta)}(l-\beta)^{\frac{1}{2}}\left(\int_{\beta}^{l}|g_{1,n}(s)|^{2}\ud s\right)^{\frac{1}{2}}\longrightarrow\,0,
\end{split}
\end{equation*}
and
\begin{equation*}
\begin{split}
\frac{1}{\fo}\left|\int_{\beta}^{l}\e^{-\fo(s-\beta)}\lambda_{n}f_{1,n}(s)\,\ud s\right|\leq C\frac{|\lambda_{n}|}{\fo^{2}}\left(f_{1,n}(\beta)-\e^{-\fo(l-\beta)}f_{1,n}(l)+\int_{\beta}^{l}|f_{1,n}'(s)|\,\ud s\right)
\\
\leq C\frac{|\lambda_{n}|}{\fo^{2}}\left(1+l-\beta+\e^{-\fo(l-\beta)}\right)\left(\max_{[\beta,l]}|f_{1,n}(s)|+\max_{[\beta,l]}|f_{1,n}'(s)|\right)\,\longrightarrow\,0,
\end{split}
\end{equation*}
\begin{equation*}
\begin{split}
\left|\int_{\beta}^{l}\e^{-\fo(2l-\beta-s)}\lambda_{n}f_{1,n}(s)\,\ud s\right|\leq C|\lambda_{n}|\e^{-\fo(l-\beta)}\max_{[\beta,l]}|f_{1,n}(s)|\,\longrightarrow\,0,
\\
\e^{-\fo(l-\beta)}\left|\int_{\beta}^{l}\cos(\fo(l-s))\lambda_{n}f_{1,n}(s)\,\ud s\right|\leq C|\lambda_{n}|\e^{-\fo(l-\beta)}\max_{[\beta,l]}|f_{1,n}(s)|\,\longrightarrow\,0.
\end{split}
\end{equation*}
We multiply the relation~\eqref{kv49} by $\e^{-\fo(x-\beta)}$ and we take $x=l$ and we use the transmission conditions $w_{1,n}(l)=w_{2,n}(l)$ and $w_{1,n}'(l)=w_{2,n}'(l)$, we have
\begin{equation}\label{kv91}
\begin{split}
\e^{-\fo(l-\beta)}w_{2,n}'(l)=\e^{-\fo(l-\beta)}w_{1,n}'(l)=-\fo\e^{-\fo(l-\beta)}w_{2,n}(l)+C_{3}
\\
+\frac{C_{2}(1-i)}{2\fo}\left[1-\e^{-\fo(1+i)(l-\beta)}\right]+\frac{1}{\fo}\int_{\beta}^{l}\int_{\beta}^{\tau}\e^{\fo(2\beta-l-\tau)}\sin(\fo(\tau-s))F_{1,n}(s)\,\ud s\,\ud\tau
\\
+\frac{C_{1}}{2\fo^{2}}\e^{-\fo(l-\beta)}\left[1-\cos(\fo(l-\beta))-\sin(\fo(l-\beta))\right].
\end{split}
\end{equation}
Since,
$$
\left|\int_{\beta}^{l}\int_{\beta}^{\tau}\e^{\fo(2\beta-l-\tau)}\sin(\fo(\tau-s))g_{1,n}(s)\,\ud s\,\ud\tau\right|\leq(l-\beta)^{\frac{3}{2}}\left(\int_{\beta	}^{l}|g_{1,n}(s)|^{2}\ud s\right)^{\frac{1}{2}}\longrightarrow\,0,
$$
and
$$
\left|\int_{\beta}^{l}\!\!\!\int_{\beta}^{\tau}\!\!\e^{\fo(2\beta-l-\tau)}\sin(\fo(\tau-s))\lambda_{n}f_{1,n}(s)\,\ud s\,\ud\tau\right|\leq C|\lambda_{n}|\e^{-\fo(l-\beta)}\max_{[\beta,l]}|f_{1,n}(s)|\longrightarrow\,0,
$$
then the fourth right hand side term of~\eqref{kv91} satisfy
\begin{equation}\label{kv106}
\int_{\beta}^{l}\int_{\beta}^{\tau}\e^{\fo(2\beta-l-\tau)}\sin(\fo(\tau-s))F_{1,n}(s)\,\ud s\,\ud\tau\,\longrightarrow\,0.
\end{equation}
The expression~\eqref{kv47} multipied by $\e^{-\fo(x-\beta)}$ and taken at the point $x=l$ leads by the transmission conditions $w_{1,n}(l)=w_{2,n}(l)$ and $q_{1}(l)w_{1,n}''(l)=q_{2}(l)w_{2,n}''(l)$ to
\begin{equation}\label{kv97}
\begin{split}
\e^{-\fo(l-\beta)}w_{2,n}''(l)=\e^{-\fo(l-\beta)}\frac{q_{1}(l)}{q_{2}(l)}w_{1,n}''(l)=\frac{q_{1}(l)}{q_{2}(l)}\fo^{2}\e^{-\fo(l-\beta)}w_{2,n}(l)
\\
+\frac{q_{1}(l)}{q_{2}(l)}C_{2}\e^{-\fo(1+i)(l-\beta)}+\frac{q_{1}(l)}{q_{2}(l)\fo}C_{1}\e^{-\fo(l-\beta)}\sin(\fo(l-\beta))
\\
+\frac{q_{1}(l)}{q_{2}(l)\fo}\e^{-\fo(l-\beta)}\int_{\beta}^{l}\sin(\fo(l-s))F_{1,n}(s)\,\ud s.
\end{split}
\end{equation}
where as done in~\eqref{kv110}, the last right hand side term of~\eqref{kv97} verify
\begin{equation}\label{kv105}
\e^{-\fo(l-\beta)}\left|\int_{\beta}^{l}\sin(\fo(l-s))F_{1,n}(s)\,\ud s\right|\,\longrightarrow\,0.
\end{equation}
Using the transmission conditions  $w_{1,n}'(l)=w_{2,n}'(l)$ and $q_{1}(l)w_{1,n}''(l)=q_{2}(l)w_{2,n}''(l)$ to substitute~\eqref{kv91} and~\eqref{kv97} into~\eqref{kv45} which multiplied  by $\e^{-\fo(x-\beta)}$ and taken in the point $x=l$, we find using $w_{1,n}(l)=w_{2,n}(l)$ and $q_{1}(l)w_{1,n}'''(l)=q_{2}(l)w_{2,n}'''(l)$ that
\begin{equation}\label{kv98}
\begin{split}
\e^{-\fo(l-\beta)}w_{2,n}'''(l)=\e^{-\fo(l-\beta)}\frac{q_{1}(l)}{q_{2}(l)}w_{1,n}'''(l)=-\frac{q_{1}(l)}{q_{2}(l)}\fo^{3}\e^{-\fo(l-\beta)}w_{2,n}(l)
\\
+\frac{q_{1}(l)}{2q_{2}(l)}\fo C_{2}\left[(1-i)-(1+i)\e^{-\fo(1+i)(l-\beta)}\right]+\frac{q_{1}(l)}{q_{2}(l)}\fo^{2}C_{3}
\\
+\frac{q_{1}(l)}{2q_{2}(l)}\e^{-\fo(l-\beta)}C_{1}\left[1-\cos(\fo(l-\beta))-(1+2i)\sin(\fo(l-\beta))+2\e^{i\fo(l-\beta)}\right]
\\
+\frac{q_{1}(l)}{q_{2}(l)}\fo\int_{\beta}^{l}\int_{\beta}^{\tau}\e^{\fo(2\beta-l-\tau)}\sin(\fo(\tau-s))F_{1,n}(s)\,\ud s\,\ud\tau
\\
-i\frac{q_{1}(l)}{q_{2}(l)}\e^{-\fo(l-\beta)}\int_{\beta}^{l}\sin(\fo(l-s))F_{1,n}(s)\,\ud s+\frac{q_{1}(l)}{q_{2}(l)}\e^{-\fo(l-\beta)}\int_{\beta}^{l}\e^{i\fo(l-s)}F_{1,n}(s)\,\ud s.
\end{split}
\end{equation}
where similarly to~\eqref{kv110}, the last term of~\eqref{kv98} verify
\begin{equation}\label{kv104}
\e^{-\fo(l-\beta)}\left|\int_{\beta}^{l}\e^{i\fo(l-s)}F_{1,n}(s)\,\ud s\right|\,\longrightarrow\,0.
\end{equation}
We substitute~\eqref{kv91},~\eqref{kv97} and~\eqref{kv98} into~\eqref{kv55} multiplied by $\e^{-\fo(l-\beta)}$, we obtain
\begin{equation}\label{kv99}
\begin{split}
K_{1}\e^{-\fo(l-\beta)}=\left[-\frac{q_{1}(l)}{q_{2}(l)}\fo^{3}+\fo\ft^{2}+i\left(\frac{q_{1}(l)}{q_{2}(l)}\fo^{2}\ft-\ft^{3}\right)\right]w_{2}(l)\e^{-\fo(l-\beta)}
\\
+C_{3}\left[\frac{q_{1}(l)}{q_{2}(l)}\fo^{2}-\ft^{2}\right]+\frac{C_{2}}{2}\bigg[(1-i)\left(\frac{q_{1}(l)}{q_{2}(l)}\fo-\fto\ft\right)
\\
+\e^{-\fo(1+i)(l-\beta)}\left(2i\frac{q_{1}(l)}{q_{2}(l)}\ft-(1+i)\frac{q_{1}(l)}{q_{2}(l)}\fo+(1-i)\fto\ft\right)\bigg]
\\
+\frac{C_{1}}{2}\e^{-\fo(l-\beta)}\bigg[\frac{2q_{1}(l)}{q_{2}(l)}\e^{i\fo(l-\beta)}-\left(\frac{q_{1}(l)}{q_{2}(l)}+\fto^{2}\right)\cos(\fo(l-\beta))
\\
+\left(\fto^{2}+2i\frac{q_{1}(l)}{q_{2}(l)}\fto-\frac{q_{1}(l)}{q_{2}(l)}(1+2i)\right)\sin(\fo(l-\beta))\bigg]
\\
+\left(\frac{q_{1}(l)}{q_{2}(l)}\fo-\fto\ft\right)\int_{\beta}^{l}\int_{\beta}^{\tau}\e^{\fo(2\beta-l-\tau)}\sin(\fo(\tau-s))F_{1,n}(s)\,\ud s\,\ud\tau
\\
+i(\fto-1)\frac{q_{1}(l)}{q_{2}(l)}\e^{-\fo(l-\beta)}\int_{\beta}^{l}\sin(\fo(l-s))F_{1,n}(s)\,\ud s
\\
+\frac{q_{1}(l)}{q_{2}(l)}\e^{-\fo(l-\beta)}\int_{\beta}^{l}\e^{i\fo(l-s)}F_{1,n}(s)\,\ud s.
\end{split}
\end{equation}
We substitute~\eqref{kv97} into~\eqref{kv56} multiplied by $\e^{-\fo(l-\beta)}$ then we obtain
\begin{equation}\label{kv100}
\begin{split}
K_{2}\e^{-\fo(l-\beta)}=\left[\frac{q_{1}(l)}{q_{2}(l)}\fo^{2}-\ft^{2}\right]w_{2}(l)\e^{-\fo(l-\beta)}+C_{2}\frac{q_{1}(l)}{q_{2}(l)}\e^{-\fo(1+i)(l-\beta)}
\\
+C_{1}\frac{q_{1}(l)}{q_{2}(l)\fo}\e^{-\fo(l-\beta)}\sin(\fo(l-\beta))+\frac{q_{1}(l)}{q_{2}(l)\fo}\e^{-\fo(l-\beta)}\int_{\beta}^{l}\sin(\fo(l-s))F_{1,n}(s)\,\ud s.
\end{split}
\end{equation}
Equation~\eqref{kv91} can be written as follow
\begin{equation}\label{kv101}
\begin{split}
(w_{2,n}'(l)+\ft w_{2,n}(l))\e^{-\fo(l-\beta)}=\left[\ft-\fo\right]w_{2,n}(l)\e^{-\fo(l-\beta)}
\\
+C_{3}+\frac{C_{1}}{2\fo^{2}}\e^{-\fo(l-\beta)}\left(1-\cos(\fo(l-\beta))-\sin(\fo(l-\beta))\right)
\\
+\frac{(1-i)}{2\fo}C_{2}\left[1-\e^{-\fo(1+i)(l-\beta)}\right]+\frac{1}{\fo}\int_{\beta}^{l}\int_{\beta}^{\tau}\e^{\fo(2\beta-l-\tau)}\sin(\fo(\tau-s))F_{1,n}(s)\,\ud s\,\ud\tau.
\end{split}
\end{equation}
Multiplying~\eqref{kv50} by $2\ft\e^{-\fo(l-\beta)}$ and taking $x=l$, we get
\begin{equation}\label{kv102}
\begin{split}
2(\ft w_{2,n}'(l)+\ft^{2}w_{2,n}(l))\e^{-\fo(l-\beta)}=K_{2}(1-i)\e^{-\fo(l-\beta)}\left[\e^{-\ft(1+i)(L-l)}-1\right]
\\
+\frac{K_{1}}{\ft}\e^{-\fo(l-\beta)}\left[-1+\e^{-\ft(L-l)}(\cos(\ft(L-l))+\sin(\ft(L-l)))\right]
\\
+2\,\e^{-\fo(l-\beta)}\int_{L}^{l}\int_{l}^{\tau}\e^{\ft(l-\tau)}\sin(\ft(\tau-s))F_{2,n}(s)\,\ud s\,\ud\tau,
\end{split}
\end{equation}
where the last term verify for the same arguments as previously that
\begin{equation}\label{kv58}
\e^{-\fo(l-\beta)}\left|\int_{L}^{l}\int_{l}^{\tau}\e^{\ft(l-\tau)}\sin(\ft(\tau-s))F_{2,n}(s)\,\ud s\,\ud\tau\right|\longrightarrow\,0.
\end{equation}
We substitute~\eqref{kv60},~\eqref{kv99},~\eqref{kv100} and~\eqref{kv101} into the relation~\eqref{kv102} then a straightforward calculations leads to
\begin{equation}\label{kv103}
\begin{split}
\frac{(1-i)}{4}\left[\frac{q_{1}(l)}{q_{2}(l)}(1+\fot)+\fto(1+\fto)\right]w_{1,n}''(\beta^{+})=O(1)\fo w_{1,n}(\beta)
\\
+O(1)\e^{-\fo(l-\beta)}\left(C_{1}+C_{2}+\e^{-\fo(l-\beta)}\omega_{n}^{2}C_{4}\right)+O(1)\omega_{n}C_{3}+O(1)\e^{-\fo(l-\beta)}
\\
+O(1)\e^{-\fo(l-\beta)}\left(\frac{1}{\ft}\int_{\beta}^{l}\e^{i\fo(l-s)}F_{1,n}(s)\,\ud s+\omega_{n}^{-1}\int_{\beta}^{l}\sin(\fo(l-s))F_{1,n}(s)\,\ud s\right)
\\
+O(1)\int_{\beta}^{l}\int_{\beta}^{\tau}\e^{\fo(2\beta-l-\tau)}\sin(\fo(\tau-s))F_{1,n}(s)\,\ud s\,\ud\tau
\\
+2\,\e^{-\fo(l-\beta)}\int_{L}^{l}\int_{l}^{\tau}\e^{\ft(l-\tau)}\sin(\ft(\tau-s))F_{2,n}(s)\,\ud s\,\ud\tau
\\
+O(1)\frac{1}{\fo}\bigg(\int_{\beta}^{l}\e^{-\fo(s-\beta)}F_{1,n}(s)\,\ud s+\int_{\beta}^{l}\e^{-\fo(2l-\beta-s)}F_{1,n}(s)\,\ud s
\\
-2\,\e^{-\fo(l-\beta)}\int_{\beta}^{l}\cos(\fo(l-s))F_{1,n}(s)\,\ud s\bigg).
\end{split}
\end{equation}
Using~\eqref{kv108},~\eqref{kv109},~\eqref{kv110},~\eqref{kv106},~\eqref{kv105},~\eqref{kv104} and~\eqref{kv58} and substituting the expressions of $C_{1},\,C_{2},\,C_{3}$ and $C_{4}$ respectively in~\eqref{kv54},~\eqref{kv92},~\eqref{kv57} and~\eqref{kv93} into~\eqref{kv103} it follows that
\begin{equation}\label{kv114}
(1+o(1))w_{1,n}''(\beta^{+})=(\omega_{n}^{2}w_{1,n}(\beta)+\omega_{n}^{2}w_{1,n}'(\beta)+\omega_{n}^{-1}w_{1,n}'''(\beta^{+})+1)o(1).
\end{equation}
From Lemma~\ref{kv111} the right hand side of~\eqref{kv114} converge to zero therefore $w_{1,n}''(\beta^{+})\,\longrightarrow\,0$. Thus we have proved the second identity of~\eqref{kv43}. Hence we proved the promised contradiction.
\end{pr}
\section{Longitudinal motion}\label{kv95}
Let $H$ be as defined in the previous section and we define $V_{2}$ by
$$V_{2}=\{u=(u_{1},u_{2})\in H^{1}(0,l)\times H^{1}(l,L):u_{1}(l)=u_{2}(l),\,u_{1}(0)=0,\,u_{2}(L)=0\}$$
with the norm
$$\|u\|_{V_{2}}=\|(u_{1},u_{2})\|_{V_{2}}=\left(\int_{0}^{l}p_{1}.|u_{1}'|^{2}(x)\,\ud x+\int_{l}^{L}p_{2}.|u_{2}'|^{2}(x)\,\ud x\right)^{\frac{1}{2}}.$$
We define $\mathcal{H}_{2}=H\times V_{2}$ with norm $\|(u,v)\|_{\mathcal{H}_{2}}^{2}=\|u\|_{V_{2}}^{2}+\|v\|_{H}^{2}$. Then $\mathcal{H}_{2}$ is a Hilbert space in which we define the operator $\mathcal{A}_{2}$ by
$$\mathcal{A}_{2}(u,v)=(v_{1},v_{2},(p_{1}u'_{1}+D_{a}v_{1}')',(p_{2}u'_{2})')$$
with domain
\begin{equation*}
\begin{split}
\mathcal{D}(\mathcal{A}_{2})=\{(u,v)=((u_{1},u_{2}),(v_{1},v_{2}))\in\mathcal{H}_{2}: v\in V_{2},\,p_{1}u_{1}'+D_{a}v_{1}\in H^{1}(0,l),\,u_{2}'\in H^{1}(l,L),
\\
p_{1}(l).u_{1}'(l)=p_{2}(l).u_{2}'(l)\}.
\end{split}
\end{equation*}
Then~\eqref{kv1} can be rewritten as an abstract evolution equation on $\mathcal{H}_{2}$,
$$(\dot{u}(t),\dot{v}(t))=\mathcal{A}_{2}(u(t),v(t)),\quad (u_{1}(0),u_{2}(0),v_{1}(0),v_{2}(0))=(u_{1}^{0},u_{2}^{0},u_{1}^{1},u_{2}^{1}).$$
By the same way as Proposition~\ref{kv96}  we can prove the following
\begin{pro}
The linear operator $\mathcal{A}_{2}$ generates a $C_{0}$-semigroup of contractions $e^{t\mathcal{A}_{2}}$ on $\mathcal{H}_{2}$, in particular there exists a unique solution of~\eqref{kv1} which can be expressed by means of a semigroup on $\mathcal{H}_{2}$ having the following regularity of the solution
$$
\left(\begin{array}{l}
u_{1}
\\
u_{2}
\\
\dot{u}_{1}
\\
\dot{u}_{2}
\end{array}\right)
\in C([0,+\infty[,\mathcal{D}(\mathcal{A}_{2}))\cap C^{1}([0,+\infty[,\mathcal{H}_{2}),$$
if $(u_{1}^{0},u_{2}^{0},u_{1}^{1},u_{2}^{1})\in \mathcal{D}(\mathcal{A}_{2})$ and a mild solution
$$
\left(\begin{array}{l}
u_{1}
\\
u_{2}
\\
\dot{u}_{1}
\\
\dot{u}_{2}
\end{array}\right)
\in C([0,+\infty[,\mathcal{H}_{2}),$$
if $(u_{1}^{0},u_{2}^{0},u_{1}^{1},u_{2}^{1})\in\mathcal{H}_{2}$. In addition the imaginary axis $i\R$ is a subset of the resolvent set $\rho(\mathcal{A}_{2})$.
\end{pro}

Assume that $p_{1}$ and $p_{2}$ are strictly positive real functions such that $p_{1}$ is constant in each of the interval $[0,\alpha)$ and $(\beta,l]$ and $p_{2}$ is constant in $(l,L]$. We suppose also that $a\in\mathrm{C}(\alpha,\beta)$.
\begin{thm}\label{kv115}
Under the above assumptions on the coefficients of~\eqref{kv1} the semigroup $e^{\mathcal{A}_{2}t}$ is polynomially stable and in particular, there exist $M>0$ such that
\[
\|e^{\mathcal{A}_{2}t}(u_{1}^{0},u_{2}^{0},u_{1}^{1},u_{2}^{1})\|_{\mathcal{H}_{2}}\leq\frac{M}{(1+t)^{2}}\|(u_{1}^{0},u_{2}^{0},u_{1}^{1},u_{2}^{1})\|_{\mathcal{D}(\mathcal{A}_{2})},\quad\forall\,(u_{1}^{0},u_{2}^{0},u_{1}^{1},u_{2}^{1})\in\mathcal{D}(\mathcal{A}_{2}),\,\forall\,t\geq 0.
\]
\end{thm}
\begin{pr}
We need only to verify the condition for a semigroup of contraction on a Hilbert space being polynomially stable (see~\cite{BT}), i.e.,
\begin{equation}\label{kv61}
\sup\{|\lambda|^{-\frac{1}{2}}.\|(i\lambda-\mathcal{A}_{2})^{-1}\|;\quad\lambda\in\R\}<+\infty.
\end{equation}
Suppose that~\eqref{kv61} is not true. By the continuity of the resolvent and the Hahn Banach theorem, there exist $\lambda_{n}\in\R_{+}$, $((u_{1,n},u_{2,n}),(v_{1,n},v_{2,n}))\in\mathcal{D}(\mathcal{A}_{2})$, $n=1,2,\ldots$, such that
\begin{equation}\label{kv62}
\|((u_{1,n},u_{2,n}),(v_{1,n},v_{2,n}))\|_{\mathcal{H}_{2}}=1,\qquad\lambda_{n}\,\longrightarrow\,+\infty
\end{equation}
and 
\begin{equation}\label{kv63}
\lambda_{n}^{\frac{1}{2}}(i\lambda_{n}-\mathcal{A}_{2})((u_{1,n},u_{2,n}),(v_{1,n},v_{2,n}))=((f_{1,n},f_{2,n}),(g_{1,n},g_{2,n}))\,\longrightarrow\,0\;\text{in}\;\mathcal{H}_{2},
\end{equation}
which mean
\begin{eqnarray}
\lambda_{n}^{\frac{1}{2}}(i\lambda_{n}u_{1,n}-v_{1,n})=f_{1,n}\,\longrightarrow\,0\;\text{in}\; H^{1}(0,l)\label{kv64}
\\
\lambda_{n}^{\frac{1}{2}}(i\lambda_{n}u_{2,n}-v_{2,n})=f_{2,n}\,\longrightarrow\,0\;\text{in}\; H^{1}(l,L)\label{kv65}
\\
\lambda_{n}^{\frac{1}{2}}(i\lambda_{n}v_{1,n}-T_{n}')=g_{1,n}\,\longrightarrow\,0\;\text{in}\; L^{2}(0,l)\label{kv66}
\\
\lambda_{n}^{\frac{1}{2}}(i\lambda_{n}v_{2,n}-(p_{2}u_{2,n})')=g_{2,n}\,\longrightarrow\,0\;\text{in}\; L^{2}(l,L)\label{kv67}
\end{eqnarray}
where
$$
T_{n}=p_{1}u_{1,n}'+D_{a}v_{1,n}'.$$
Our goal is to find a contradiction with~\eqref{kv62}.

We first consider~\eqref{kv64} and~\eqref{kv66} on the interval $(\alpha,\beta)$. From~\eqref{kv63}, we obtain
\begin{equation}\label{kv68}
\lambda_{n}^{\frac{1}{2}}\int_{\alpha}^{\beta}\!\!a|v_{1,n}'|^{2}\ud x=\re\left\langle\lambda_{n}^{\frac{1}{2}}(i\lambda_{n}-\mathcal{A}_{2})(u_{1,n},u_{2,n},v_{1,n},v_{2,n}),(u_{1,n},u_{2,n},v_{1,n},v_{2,n}) \right\rangle_{\mathcal{H}_{2}}\longrightarrow0
\end{equation}
which imply that
\begin{equation}\label{kv69}
\|\lambda_{n}^{\frac{1}{4}}v_{1,n}'\|_{L^{2}(\alpha,\beta)}\,\longrightarrow\,0,\quad\text{and}\quad \|\lambda_{n}^{\frac{5}{4}}u_{1,n}'\|_{L^{2}(\alpha,\beta)}\,\longrightarrow\,0.
\end{equation}
Thus we also have
\begin{equation}\label{kv70}
\|\lambda_{n}^{\frac{1}{4}}T_{n}\|_{L^{2}(\alpha,\beta)}\,\longrightarrow\,0.
\end{equation}
The rest of the proof depend on the following lemma.
\begin{lem}
The functions $u_{1,n}$ and $T_{n}$ have the following properties:
\begin{eqnarray}
&|T_{n}(\alpha^{+})|\,\longrightarrow\,0,&\qquad |T_{n}(\beta^{-})|\,\longrightarrow\,0,\label{kv71}
\\
&|\lambda_{n}u_{1,n}(\alpha^{+})|\,\longrightarrow\,0,&\qquad |\lambda_{n}u_{1,n}(\beta^{-})|\,\longrightarrow\,0.\label{kv72}
\end{eqnarray}
\end{lem}
\begin{pr}
From~\eqref{kv64} we have
\begin{equation}\label{kv73}
\lambda_{n}^{-1}\|\psi v_{1,n}\|_{H^{1}(\alpha,\beta)}=O(1),\qquad\forall\,\psi\in\mathrm{C}^{\infty}(0,l),\;\supp(\psi)\subset(\alpha,\beta).
\end{equation}
Equations~\eqref{kv66},~\eqref{kv70} and~\eqref{kv73} imply that
\begin{equation}\label{kv74}
\lambda_{n}^{\frac{1}{2}}\int_{\alpha}^{\beta}\psi|v_{1,n}|^{2}\,\ud x\,\longrightarrow\,0,\qquad\forall\,\psi\in\mathrm{C}^{\infty}(0,l),\;\supp(\psi)\subset(\alpha,\beta).
\end{equation}
Applying the interpolation theorem involving subdomains~\cite[Theorem 4.23]{A} we find that~\eqref{kv69} and~\eqref{kv74} imply
\begin{equation}\label{kv75}
\lambda_{n}^{\frac{1}{4}}v_{1,n}\,\longrightarrow\,0\qquad\text{in }H^{1}(\alpha,\beta).
\end{equation}
We take the inner product of~\eqref{kv66} with $v_{1,n}$ in $L^{2}(\alpha,\beta)$ to obtain
\begin{equation}\label{kv76}
\begin{split}
\int_{\alpha}^{\beta}g_{1,n}.\overline{v}_{1,n}\,\ud x=i\lambda_{n}^{\frac{3}{2}}\int_{\alpha}^{\beta}|v_{1,n}|^{2}\,\ud x+\int_{\alpha}^{\beta}\lambda_{n}^{\frac{1}{4}}T_{n}.\lambda_{n}^{\frac{1}{4}}\overline{v}_{1,n}'\,\ud x
\\
+\lambda_{n}^{\frac{1}{2}}T_{n}(\alpha^{+})\overline{v}_{1,n}(\alpha^{+})-\lambda_{n}^{\frac{1}{2}}T_{n}(\beta^{-})\overline{v}_{1,n}(\beta^{-}).
\end{split}
\end{equation}
Using~\eqref{kv66},~\eqref{kv70},~\eqref{kv75} and the interpolation inequality we obtain for the third term in the right hand side of~\eqref{kv76} that
\begin{eqnarray}\label{kv77}
\lambda_{n}^{\frac{1}{2}}.|T_{n}(\alpha^{+})|.|v_{1,n}(\alpha^{+})|\!\!\!\!\!&\leq&\!\!\!C\lambda_{n}^{\frac{1}{2}}\|T_{n}\|_{L^{2}(\alpha,\beta)}^{\frac{1}{2}}.\|T_{n}'\|_{L^{2}(\alpha,\beta)}^{\frac{1}{2}}.\|v_{1,n}\|_{L^{2}(\alpha,\beta)}^{\frac{1}{2}}.\|v_{1,n}'\|_{L^{2}(\alpha,\beta)}^{\frac{1}{2}}\nonumber
\\
&=&\!\!\!C\|\lambda_{n}^{\frac{1}{4}}T_{n}\|_{L^{2}(\alpha,\beta)}^{\frac{1}{2}}.\|\lambda_{n}^{\frac{1}{4}}v_{1,n}\|_{L^{2}(\alpha,\beta)}^{\frac{1}{2}}.\|\lambda_{n}^{\frac{1}{4}}v_{1,n}'\|_{L^{2}(\alpha,\beta)}^{\frac{1}{2}}+o(1)\nonumber
\\
&=& o(1)(1+\|\lambda^{\frac{3}{4}}v_{1,n}\|_{L^{2}(\alpha,\beta)}).
\end{eqnarray}
And similarly we can show that
\begin{equation}\label{kv78}
|\lambda_{n}^{\frac{1}{2}}|.|T_{n}(\beta^{-})|.|v_{1,n}(\beta^{-})|\leq o(1)(1+\|\lambda^{\frac{3}{4}}v_{1,n}\|_{L^{2}(\alpha,\beta)}).
\end{equation}
Now~\eqref{kv76},~\eqref{kv77} and~\eqref{kv78} leads to
\begin{equation}\label{kv79}
\|\lambda_{n}^{\frac{3}{4}}v_{1,n}\|_{L^{2}(\alpha,\beta)}\,\longrightarrow\,0.
\end{equation}
We multiply~\eqref{kv66} by $\lambda_{n}^{-\frac{1}{2}}(\beta-x)\overline{T}_{n}$ and we take the inner product in $L^{2}(\alpha,\beta)$ we find
\begin{equation}\label{kv80}
\re\left\langle i\lambda_{n}^{\frac{3}{4}}v_{1,n},\lambda_{n}^{\frac{1}{4}}(\beta-x)T_{n}\right\rangle_{L^{2}(\alpha,\beta)}-\frac{1}{2}(\beta-\alpha)|T_{n}(\alpha^{+})|^{2}-\frac{1}{2}\|T_{n}\|_{L^{2}(\alpha,\beta)}^{2}\,\longrightarrow\,0,
\end{equation}
and by multiplying~\eqref{kv66} by $\lambda_{n}^{-\frac{1}{2}}(\alpha-x)\overline{T}_{n}$ and taking the inner product in $L^{2}(\alpha,\beta)$ we get
\begin{equation}\label{kv81}
\re\left\langle i\lambda_{n}^{\frac{3}{4}}v_{1,n},\lambda_{n}^{\frac{1}{4}}(\alpha-x)T_{n}\right\rangle_{L^{2}(\alpha,\beta)}-\frac{1}{2}(\beta-\alpha)|T_{n}(\beta^{-})|^{2}-\frac{1}{2}\|T_{n}\|_{L^{2}(\alpha,\beta)}^{2}\,\longrightarrow\,0.
\end{equation}
Hence~\eqref{kv71} follows now since the first and the third terms of~\eqref{kv80} and~\eqref{kv81} converge to zero by~\eqref{kv70} and~\eqref{kv79}. On the other hand by~\eqref{kv65} and~\eqref{kv75} we obtain
\begin{equation}\label{kv82}
\|\lambda_{n}^{\frac{5}{4}}u_{n}\|_{L^{2}(\alpha,\beta)}\,\longrightarrow\,0,
\end{equation}
then by~\eqref{kv69} and~\eqref{kv82} we get
$$
\|\lambda_{n}^{\frac{5}{4}}u_{n}\|_{H^{1}(\alpha,\beta)}\,\longrightarrow\,0,
$$
hence~\eqref{kv72} hold from the Sobolev embedding inequalities.
\end{pr}

Using the continuity conditions at $x=\alpha$ and $x=\beta$, we arrive form~\eqref{kv71} and~\eqref{kv72} at
\begin{equation}\label{kv83}
\begin{split}
|T_{n}(\alpha^{-})|\,\longrightarrow\,0,&\qquad |T_{n}(\beta^{+})|\,\longrightarrow\,0,
\\
|\lambda_{n}u_{1,n}(\alpha^{-})|\,\longrightarrow\,0,&\qquad |\lambda_{n}u_{1,n}(\beta^{+})|\,\longrightarrow\,0.
\end{split}
\end{equation}
We consider now~\eqref{kv64}-\eqref{kv67} on the intervals $(0,\alpha)$, $(\beta,l)$ and $(l,L)$, then by replacing~\eqref{kv64} and~\eqref{kv65} respectively into~\eqref{kv66} and~\eqref{kv67} we obtain
\begin{eqnarray}
-\lambda_{n}^{2}u_{1,n}-(p_{1}u_{1,n}')'=\lambda_{n}^{-\frac{1}{2}}g_{1,n}+i\lambda_{n}^{\frac{1}{2}}f_{1,n}&\text{in}&(0,\alpha)\cup(\beta,l),\label{kv84}
\\
-\lambda_{n}^{2}u_{2,n}-(p_{2}u_{2,n}')'=\lambda_{n}^{-\frac{1}{2}}g_{2,n}+i\lambda_{n}^{\frac{1}{2}}f_{2,n}&\text{in}&(l,L).\label{kv85}
\end{eqnarray}
Take the inner product of~\eqref{kv84} with $\kappa(x)u_{1,n}'$ in $L^{2}(0,\alpha)$ where $\kappa\in H^{1}(0,\alpha)$ and $\kappa(0)=0$, then with $\kappa_{1}(x)u_{1,n}'$ in $L^{2}(\beta,l)$ where $\kappa_{1}\in H^{1}(\beta,l)$ , and the inner product of~\eqref{kv85} with $\kappa_{2}(x)u_{2,n}'$ in $L^{2}(l,L)$ where $\kappa_{2}\in H^{1}(l,L)$ and $\kappa_{2}(L)=0$. A straightforward calculation shows that the real part of this inner products leads to the following
\begin{equation}\label{kv86}
\int_{0}^{\alpha}\kappa'.|\lambda_{n}u_{1,n}|^{2}\,\ud x+\int_{0}^{\alpha}\kappa'p_{1}.|u_{1,n}'|^{2}\,\ud x=\kappa(\alpha)|\lambda_{n}u_{1,n}(\alpha)|^{2}+\kappa(\alpha)|u_{1,n}'(\alpha^{-})|^{2}+o(1),
\end{equation}
\begin{equation}\label{kv87}
\begin{split}
\int_{\beta}^{l}\kappa_{1}'.|\lambda_{n}u_{1,n}|^{2}\,\ud x+\int_{\beta}^{l}\kappa_{1}'p_{1}.|u_{1,n}'|^{2}\,\ud x=\kappa_{1}(l)|\lambda_{n}u_{1,n}(l)|^{2}
\\
+\kappa_{1}(l)p_{1}(l)|u_{1,n}'(l)|^{2}-\kappa_{1}(\beta)|\lambda_{n}u_{1,n}(\beta)|^{2}-\kappa_{1}(\beta)p_{1}(\beta)|u_{1,n}'(\beta^{+})|^{2}+o(1),
\end{split}
\end{equation}
and
\begin{equation}\label{kv88}
\int_{l}^{L}\kappa_{2}'.|\lambda_{n}u_{2,n}|^{2}\,\ud x+\int_{l}^{L}\kappa_{2}'p_{2}.|u_{2,n}'|^{2}\,\ud x=-\kappa_{2}(l)|\lambda_{n}u_{2,n}(l)|^{2}-\kappa_{2}(l)p_{2}(l)|u_{2,n}'(l)|^{2}+o(1).
\end{equation}
Moreover we can let
\begin{eqnarray*}
\kappa(x)=x\qquad\kappa_{1}(x)=x-\gamma\quad\text{and}\quad\kappa_{2}(x)=x-L
\end{eqnarray*}
that's verify
$$
\kappa(0)=\kappa_{2}(L)=0
$$
and where $\gamma$ is chosen such that 
\begin{equation}\label{kv89}
\gamma\geq L\qquad\text{and}\qquad\gamma\geq l+(L-l)\frac{p_{1}(l)}{p_{2}(l)} 
\end{equation}
Summing~\eqref{kv86}-\eqref{kv88} and using~\eqref{kv83} and~\eqref{kv89} then we obtain by transmission conditions
\begin{equation*}
\begin{split}
\int_{0}^{\alpha}|\lambda_{n}u_{1,n}|^{2}\,\ud x+\int_{\beta}^{l}|\lambda_{n}u_{1,n}|^{2}\,\ud x+\int_{l}^{L}|\lambda_{n}u_{2,n}|^{2}\,\ud x+\int_{0}^{\alpha}p_{1}|u_{1,n}'|^{2}\,\ud x
\\
+\int_{\beta}^{l}p_{1}|u_{1,n}'|^{2}\,\ud x+\int_{l}^{L}p_{2}|u_{2,n}'|^{2}\,\ud x=(L-\gamma)|\lambda_{n}u_{n}|^{2}(l)
\\
+p_{1}(l)\left(l-\gamma+(L-l)\frac{p_{1}(l)}{p_{2}(l)}\right)|u_{1,n}'|^{2}(l)+o(1).
\end{split}
\end{equation*}
Then by~\eqref{kv89} we show that 
\begin{equation*}
\begin{split}
\int_{0}^{\alpha}|\lambda_{n}u_{1,n}|^{2}\,\ud x+\int_{\beta}^{l}|\lambda_{n}u_{1,n}|^{2}\,\ud x+\int_{l}^{L}|\lambda_{n}u_{2,n}|^{2}\,\ud x+
\\
\int_{0}^{\alpha}p_{1}|u_{1,n}'|^{2}\,\ud x+\int_{\beta}^{l}p_{1}|u_{1,n}'|^{2}\,\ud x+\int_{l}^{L}p_{2}|u_{2,n}'|^{2}\,\ud x\,\longrightarrow\,0.
\end{split}
\end{equation*}
Hence, it follows from~\eqref{kv64} and~\eqref{kv65} that
\begin{equation}\label{kv90}
\begin{split}
\int_{0}^{\alpha}|v_{1,n}|^{2}\,\ud x+\int_{\beta}^{l}|v_{1,n}|^{2}\,\ud x+\int_{l}^{L}|v_{2,n}|^{2}\,\ud x+
\\
\int_{0}^{\alpha}p_{1}|u_{1,n}'|^{2}\,\ud x+\int_{\beta}^{l}p_{1}|u_{1,n}'|^{2}\,\ud x+\int_{l}^{L}p_{2}|u_{2,n}'|^{2}\,\ud x\,\longrightarrow\,0.
\end{split}
\end{equation}
Finally, we combine~\eqref{kv90},~\eqref{kv75} and~\eqref{kv69} we obtain the promised contradiction.
\end{pr}

Next we show that the polynomial decay rate given in Theorem~\ref{kv115} is sharp.
The main idea of the proof is to show that the resolvent $\displaystyle\lambda^{\epsilon-\frac{1}{2}}(i\lambda-\mathcal{A}_{2})^{-1}$ is not uniformly bounded with respect to $\lambda\in\R$ where $\displaystyle 0<\epsilon\leq\frac{1}{2}$.
Let $\displaystyle\lambda=\lambda_{n}=\frac{2n\pi\sqrt{p_{2}}}{(L-l)}$, $n=1,2,\dots,$. We take $\displaystyle\alpha=l-\beta=\frac{(L-l)\sqrt{p_{1}}}{\sqrt{p_{2}}}$, $p_{1}$, $p_{2}$ and $a$ constants in $(0,l)$, $(l,L)$ and $(\alpha,\beta)$ respectively and we define
\begin{equation*}
f(x)=f_{n}(x)=\left\{\begin{array}{ll}
0&0<x<l
\\
\displaystyle\frac{1}{\lambda}\sin\left(\frac{\lambda(x-l)}{\sqrt{p_{2}}}\right)&l<x<L
\end{array}\right.
\end{equation*}
and
\begin{equation*}
g(x)=g_{n}(x)=\left\{\begin{array}{ll}
0&0<x<l
\\
\displaystyle\cos\left(\frac{\lambda(x-l)}{\sqrt{p_{2}}}\right)&l<x<L.
\end{array}\right.
\end{equation*}
We set $f_{1,n},\,f_{2,n}$ the restriction of $f_{n}$ over the intervals $(0,l)$ and $(l,L)$ respectively and $g_{1,n},\,g_{2,n}$ the restriction of $g_{n}$ over the intervals $(0,l)$ and $(l,L)$ respectively. In the intervals $(0,\alpha),\,(\alpha,\beta),\,(\beta,l)$ and $(l,L)$ we solve the resolvent equation
$$
\lambda^{\frac{1}{2}-\epsilon}(i\lambda_{n}-\mathcal{A}_{2})\left(\begin{array}{l}
u_{1,n}
\\
u_{2,n}
\\
v_{1,n}
\\
v_{2,n}
\end{array}\right)=\left(\begin{array}{l}
f_{1,n}
\\
f_{2,n}
\\
g_{1,n}
\\
g_{2,n}
\end{array}\right)
$$
where $(u,v)=(u_{n},v_{n})=((u_{1,n},u_{2,n}),(v_{1,n},v_{2,n}))\in\mathcal{D}(\mathcal{A}_{2})$.

For $x\in(0,\alpha)$, we have
\begin{equation*}
\left\{\begin{array}{l}
i\lambda u-v=0
\\
i\lambda v-p_{1}u''=0
\\
u(0)=0,
\end{array}\right.
\end{equation*}
where the solution is given by
$$
u(x)=c_{1}\sin\left(\frac{\lambda x}{\sqrt{p_{1}}}\right)
$$

For $x\in(\alpha,\beta)$, we have
\begin{equation*}
\left\{\begin{array}{l}
i\lambda u-v=0
\\
i\lambda v-(p_{1}+ia\lambda)u''=0,
\end{array}\right.
\end{equation*}
where the solution is given by
$$
u(x)=c_{2}\e^{\omega x}+c_{3}\e^{-\omega x}
$$
with
$$
\omega=\omega_{n}=\frac{\lambda}{(p_{1}^{2}+a^{2}\lambda^{2})^{\frac{1}{4}}}\left(\cos\left(\frac{\theta}{2}\right)+i\sin\left(\frac{\theta}{2}\right)\right)
$$
where
$$
\cos(\theta)=\frac{-p_{1}}{(p_{1}^{2}+a^{2}\lambda^{2})^{\frac{1}{2}}}\longrightarrow0\;\text{ and }\;\sin(\theta)=\frac{a\lambda}{(p_{1}^{2}+a^{2}\lambda^{2})^{\frac{1}{2}}}\longrightarrow1\text{ as }n\longrightarrow+\infty.
$$
By continuity condition at $x=\alpha$,
\begin{equation*}
\left\{\begin{array}{l}
u(\alpha^{-})=u(\alpha^{-})
\\
p_{1}u'(\alpha^{-})=(p_{1}+ia\lambda)u'(\alpha^{+})
\end{array}\right.
\end{equation*}
and taking into account the expression of $\alpha$ we can solve $c_{2}$ and $c_{3}$ to get
$$
u(x)=\frac{c_{1}\lambda\sqrt{p_{1}}}{\omega(p_{1}+ia\lambda)}\sinh(\omega(x-\alpha)).
$$
Therefore,
$$
u(\beta^{-})=\frac{c_{1}\lambda\sqrt{p_{1}}}{\omega(p_{1}+ia\lambda)}\sinh(\omega(\beta-\alpha))\;\text{ and }\;u'(\beta^{-})=\frac{c_{1}\lambda\sqrt{p_{1}}}{(p_{1}+ia\lambda)}\cosh(\omega(\beta-\alpha)).
$$

For $x\in(\beta,l)$, we have
\begin{equation*}
\left\{\begin{array}{l}
i\lambda u-v=0
\\
i\lambda v-p_{1}u''=0,
\end{array}\right.
\end{equation*}
where the solution is given by
$$
u(x)=c_{4}\cos\left(\frac{\lambda x}{\sqrt{p_{1}}}\right)+c_{5}\sin\left(\frac{\lambda x}{\sqrt{p_{1}}}\right).
$$
By continuity conditions at $x=\beta$
\begin{equation*}
\left\{\begin{array}{l}
u(\beta^{-})=u(\beta^{+})
\\
(p_{1}+ia\lambda)u'(\beta^{-})=u'(\beta^{+})
\end{array}\right.
\end{equation*}
we can solve $c_{4}$ and $c_{5}$ to get
$$
u(x)=c_{1}\left[\frac{\lambda\sqrt{p_{1}}}{\omega(p_{1}+ia\lambda)}\cos\left(\frac{\lambda(x-\beta)}{\sqrt{p_{1}}}\right)\sinh(\omega(\beta-\alpha))+p_{1}\sin\left(\frac{\lambda(x-\beta)}{\sqrt{p_{1}}}\right)\cosh(\omega(\beta-\alpha))\right].
$$
Therefore, by using the expression of $l-\beta$ gived above we get
$$
u(l^{-})=\frac{c_{1}\lambda\sqrt{p_{1}}}{\omega(p_{1}+ia\lambda)}\sinh(\omega(\beta-\alpha))\;\text{ and }\;u'(l^{-})=c_{1}\lambda\sqrt{p_{1}}\cosh(\omega(\beta-\alpha)).
$$

For $x\in(l,L)$, we have
\begin{equation}\label{kv112}
\left\{\begin{array}{l}
i\lambda u-v=\lambda^{\epsilon-\frac{1}{2}}f
\\
i\lambda v-p_{2}u''=\lambda^{\epsilon-\frac{1}{2}}g
\\
u(L)=0.
\end{array}\right.
\end{equation}
Let $\displaystyle z_{\pm}=\frac{1}{2}(v(x)\pm\sqrt{p_{2}} u'(x))$. Then~\eqref{kv112} can be transformed into the first-order, diagonal and non homogeneous system in $(l,L)$ as follow
$$
\left(\begin{array}{l}
z_{+}
\\
z_{-}
\end{array}\right)'=\left(\begin{array}{lr}
\frac{i\lambda}{\sqrt{p_{2}}}&0
\\
0&-\frac{i\lambda}{\sqrt{p_{2}}}
\end{array}\right)\left(\begin{array}{l}
z_{+}
\\
z_{-}
\end{array}\right)+\left(\begin{array}{c}
-\frac{\lambda^{\epsilon-\frac{1}{2}}}{\sqrt{p_{2}}}
\\
0
\end{array}\right)g:=A\left(\begin{array}{c}
z_{+}
\\
z_{-}
\end{array}\right)+\left(\begin{array}{c}
-\frac{\lambda^{\epsilon-\frac{1}{2}}}{\sqrt{p_{2}}}
\\
0
\end{array}\right)g.
$$
Using the boundary condition $v(L)=0=z_{+}(L)+z_{-}(L)$ we obtain the solution
$$
\left(\begin{array}{l}
z_{+}
\\
z_{-}
\end{array}\right)=z_{+}(L)\e^{(x-L)A}\left(\begin{array}{r}
1
\\
-1
\end{array}\right)+\int_{x}^{L}\e^{(x-\tau)A}\left(\begin{array}{c}
\frac{\lambda^{\epsilon-\frac{1}{2}}}{\sqrt{p_{2}}}
\\
0
\end{array}\right)g(\tau)\,\ud \tau
$$
Using the fact that $v(x)=z_{+}(x)+z_{-}(x)$, we follow
\begin{equation}\label{kv113}
v(x)=2iz_{+}(L)\sin\left(\frac{2n\pi(x-L)}{L-l}\right)+\frac{\lambda^{\epsilon-\frac{1}{2}}}{\sqrt{p_{2}}}\int_{x}^{L}\e^{\frac{2in\pi(x-\tau)}{L-l}}\cos\left(\frac{2n\pi(\tau-l)}{L-l}\right)\,\ud\tau,
\end{equation}
and in particular, we have
$$
v(l^{+})=\frac{L-l}{2\sqrt{p_{2}}}\lambda^{\epsilon-\frac{1}{2}}.
$$
Furthermore, since $i\lambda u(l^{+})=f(l^{+})+v(l^{+})=v(l^{+})$ then we find that
$$
u(l^{+})=\frac{L-l}{2i\lambda\sqrt{p_{2}}}\lambda^{\epsilon-\frac{1}{2}}.
$$
Similarly, using $\sqrt{p_{2}}u'(x)=z_{+}(x)-z_{-}(x)$ then we get
$$
p_{2}u'(l^{+})=2\sqrt{p_{2}}z_{+}(L)+\frac{L-l}{2}\lambda^{\epsilon-\frac{1}{2}}.
$$
Using the transmission conditions
\begin{equation*}
\left\{\begin{array}{l}
u(l^{-})=u(l^{+})
\\
p_{1}u'(l^{-})=p_{2}u'(l^{+})
\end{array}\right.
\end{equation*}
we can solve $c_{1}$ and we find that
$$
z_{+}(L)=\lambda^{\epsilon-\frac{1}{2}}\frac{\omega p_{1}(p_{1}+ia\lambda)(L-l)}{4i\lambda p_{2}}\coth(\omega(\beta-\alpha))-\frac{L-l}{4\sqrt{p_{2}}}\lambda^{\epsilon-\frac{1}{2}}.
$$
It follows from the definition of $\lambda$ and $\omega$ that
$$
\lambda\longrightarrow+\infty,\;\;\;|\omega|\longrightarrow+\infty\;\text{ and }\;\coth(\omega(\beta-\alpha))\longrightarrow1,
$$
hence we get
$$
|z_{+}(L)|\longrightarrow+\infty\text{ as }n\longrightarrow+\infty.
$$
Referring to~\eqref{kv113} this leads to
\begin{equation*}
\begin{split}
\int_{l}^{L}|v_{2,n}(x)|^{2}\ud x&\geq 4|z_{+}(L)|^{2}\int_{l}^{L}\sin^{2}\left(\frac{2n\pi(x-L)}{L-l}\right)\ud x-\frac{(L-l)^{2}}{2p_{2}}\lambda^{2\epsilon-1}
\\
&\geq2(L-l)|z_{+}(L)|^{2}-\frac{(L-l)^{2}}{2p_{2}}\lambda^{2\epsilon-1}\longrightarrow+\infty\text{ as }n\longrightarrow+\infty.
\end{split}
\end{equation*}
Since $\|(f_{n},g_{n})\|_{\mathcal{H}_{2}}^{2}=(L-l)$ and
$$
\|\lambda^{\epsilon-\frac{1}{2}}(i\lambda_{n}-\mathcal{A}_{2})^{-1}(f_{n},g_{n})\|_{\mathcal{H}_{2}}^{2}=\|(u_{n},v_{n})\|_{\mathcal{H}_{2}}^{2}\geq\int_{l}^{L}|v_{2,n}(x)|^{2}\ud x\longrightarrow+\infty \text{ as } n\longrightarrow+\infty,
$$
we conclude that $\displaystyle\sup_{\lambda\in\R}\|\lambda^{\epsilon-\frac{1}{2}}(i\lambda-\mathcal{A}_{2})^{-1}\|_{\mathcal{L}(\mathcal{H}_{2})}=+\infty$, thus we get the optimality of the polynomial decay rate.
\subsubsection*{Acknowledgments}
The author thanks the referees for many valuable remarks which helped us to improve the paper significantly.
\nocite{*}
\bibliographystyle{alpha}
\bibliography{BibKV}
\addcontentsline{toc}{section}{References}
\end{document}